\documentclass[runningheads]{llncs}
\usepackage[russian, english]{babel}

\usepackage{graphicx}
\usepackage{float}
\usepackage{mathtools}
\usepackage{amssymb, amsmath, latexsym, bm}
\usepackage{nicefrac}
\usepackage{hhline}
\usepackage{multirow}
\usepackage{pifont}
\newcommand{\specialcell}[2][c]{%
  \begin{tabular}[#1]{@{}c@{}}#2\end{tabular}}
\usepackage[colorinlistoftodos,bordercolor=orange,backgroundcolor=orange!20,linecolor=orange,textsize=scriptsize]{todonotes}
\usepackage{algorithm}
\usepackage{algpseudocode}
\usepackage{hyperref}

\usepackage[bibliography=common]{apxproof}
\newtheoremrep{thm}{Theorem}
\newtheoremrep{prop}{Proposition}
\newtheoremrep{lmm}{Lemma}

\newcommand{\argmin}{\mathop{\arg\!\min}}
\newcommand{\argmax}{\mathop{\arg\!\max}}

\def\ee{\mathbf{e}}
\def\R{\mathbb{R}}

\newtheorem{approach}{Approach}

\begin{document}
\title{Solving smooth min-min and min-max problems by mixed oracle algorithms\thanks{The research of A. Gasnikov and P. Dvurechensky was supported by Russian Science Foundation (project No. 21-71-30005). The research of E. Gladin, A. Sadiev and A. Beznosikov was partially supported by Andrei Raigorodskii scholarship.}}
%
%
\author{Egor Gladin\inst{1,2}\orcidID{0000-0002-9086-996X} \and
Abdurakhmon Sadiev\inst{1} \and
Alexander Gasnikov\inst{1,3,5}\orcidID{0000-0002-7386-039X} \and
Pavel Dvurechensky\inst{4,3,5}\orcidID{0000-0003-1201-2343} \and
Aleksandr Beznosikov\inst{1,3}\orcidID{0000-0002-3217-3614} \and
Mohammad Alkousa\inst{1}\orcidID{0000-0001-5470-0182}}
\authorrunning{E. Gladin et al.}
%
\institute{Moscow Institute of Physics and Technology, Russia \and Skolkovo Institute of Science and Technology, Russia \and
HSE University, Russia \and
Weierstrass Institute for Applied Analysis and Stochastics, Germany
\and Institute for Information Transmission Problems RAS, Russia 
}
\maketitle              
\begin{abstract}
    In this paper we consider two types of problems which have some similarity in their structure, namely, min-min problems and min-max saddle-point problems. Our approach is based on considering the outer minimization problem as a minimization problem with inexact oracle. This inexact oracle is calculated via inexact solution of the inner problem, which is either a minimization or a maximization problem. 
Our main assumptions are that the problem is smooth and the available oracle is mixed: it is only possible to evaluate the gradient w.r.t. the outer block of variables which corresponds to the outer minimization problem, whereas for the inner problem only zeroth-order oracle is available. To solve the inner problem we use accelerated gradient-free method with zeroth-order oracle. To solve the outer problem we use either inexact variant of the Vaydya's cutting-plane method or a variant of accelerated gradient method. As a result we propose a framework which leads to non-asymptotic complexity bounds for both min-min and min-max problems. Moreover, we estimate separately the number of first- and zeroth-order oracle calls which are sufficient to reach any desired accuracy.

\end{abstract}
\section{Introduction}
In this paper, we consider smooth optimization problems in which the decision variable is decomposed into two blocks with minimization w.r.t. one block, which we call the outer block, and two types of operations w.r.t. the second block, which we call the inner block: minimization or maximization. In other words, we consider smooth min-min problems and min-max problems. The main difference of our setting with existing in the literature is that we assume that it is possible to evaluate the gradient w.r.t. the outer block of the variables, i.e. first-order oracle, and only function values, i.e. zeroth-order oracle, when we deal with the inner block of variables. Thus, we operate with mixed type of oracle: first-order in one block of variables and zeroth-order in the second block of variables.

Our motivation, firstly, comes from min-max saddle-point 
problems, which have recently became of an increased interest in machine learning community in application to training Generative Adversarial Networks \cite{goodfellow2014generative}, and other adversarial models \cite{madry2018towards}, as well as to robust reinforcement learning \cite{pinto2017robust}. The standard process is to simultaneously train neural network, find adversarial examples and learn the network to distinguish the true examples from the artificially generated. In the training process the gradient is available through the backpropagation, whereas for the generating adversarial examples the network is available as a black box and only zeroth-order oracle is available. Another close application area is 
Adversarial Attacks \cite{goodfellow2014explaining,tramr2017ensemble} on neural networks, in particular the Black-Box Adversarial Attacks \cite{conf/cvpr/NarodytskaK17}. Here the goal is for a trained network  to find a perturbation of the data in such a way that the network outputs wrong prediction. Then the training is repeated to make the network robust to such attacks. 
Since the attacking model does not have access to the architecture of the main network, but only to the input and output of the network, the only available oracle for the attacker is the zeroth-order oracle for the loss function. 
The motivation for min-min problems comes from simulation optimization \cite{fu2015handbook,shashaani2018ASTRO}, where some parts of the optimized system can be given as a black box with unavailable or computationally expensive gradients, and other parts of the objective are differentiable.

Separately zeroth-order \cite{conn2009introduction} and first-order  \cite{nesterov2018lectures} are very well developed areas of modern numerical optimization. There are also plenty of works on first-order \cite{Korpelevich1976TheEM,Nemirovski2004,nedic2009subgradient,chambolle2011first-order} and zeroth-order methods \cite{wang2020zerothorder,liu2019minmax,beznosikov_sadiev_gasnikov,sadiev2020zeroth-order} for saddle-point problems. Our main idea and contribution in this paper is to consider mixed oracles, which seems to be an underdeveloped area of optimization and saddle-point problems. In \cite{sadiev2020zeroth-order} the authors consider methods with mixed oracle, but, unlike this work, only in the context of saddle-point problems and without acceleration techniques.

Notably, in this paper we develop a generic approach which is suitable for both types of problems: min-min and min-max, and is based on the same idea for both problems: we consider the minimization problem w.r.t. the outer group of variables as a minimization problem with inexact oracle. This inexact oracle is evaluated via inexact solution of the inner problem, which is either a minimization or a maximization problem. We carefully estimate with what accuracy one needs to solve the inner problem to be able to solve the outer problem with the desired accuracy. Moreover, we have to account for the random nature of the solution to the inner problem since we use randomized gradient-free methods with zeroth-order oracle to solve the inner problem. In our approach we consider two settings for the outer problem. If the dimension of the outer problem is small, we use Vaydya's cutting-plane method  \cite{vaidya1996new,vaidya1989new}, for which we extend the analysis to the case of approximate subgradients. The drawback of this method is that it scales quite badly with the dimension. Thus, if the dimension of the outer problem is large, we exploit accelerated gradient method, for which we develop an analysis in the case when an inexact oracle is available only with some probability, which may be of independent interest. Our approach based on inner-outer loops allows also to separate complexities, i.e. the number of calls to each of the oracles: first-order oracle for the outer block of variables and zeroth-order oracle for the inner  block of variables. 

The rest of the paper is organized as follows. First we consider min-min problems in two settings: small and large dimension of the outer problem. In the first case we develop an inexact variant of the  Vaydya's method use it in the outer loop in a combination with accelerated random gradient-free method in the inner loop. In the second case, we apply accelerated gradient method in the outer loop combined with the same method in the inner loop. After that we consider saddle-point min-max problems again in two settings. When the dimension of the outer problem is small, we use the same scheme with inexact Vaydya's method and accelerated random gradient-free method. The situation is more complicated when the dimension of the outer problem is large. In this case we use a three-loop structure with the Catalyst acceleration scheme \cite{catalyst} combined with accelerated gradient method with inexact oracle and accelerated random gradient-free method.

\begin{table}[H]
\caption{Main results}
\label{tab:hresult}
\centering
\begin{tabular}{ || c c || c c ||}
\hline\hline
& & \multicolumn{2}{c||}{{\large\bfseries Oracle Complexity}}\\ [2ex]
& & $\mathbf{1}$-st order & $\mathbf{0}$-th order\\ [1ex]
\hline\hline 
\multirow{1}{*}{{\bfseries \specialcell{Min \\Min}}}& \specialcell{Small \\Scale} & $\widetilde{O}\left(n_x\right)$ & $\widetilde{O}\left(n_x n_{y} \sqrt{\frac{L_{y y}}{\mu_{y}}}\right)$ \\ [4ex] 
\hline \hline
\multirow{2}{*}{{\bfseries \specialcell{Min \\Max}}}& \specialcell{Small \\Scale} & $\widetilde{O}\left(n_x\right)$ & $\widetilde{O}\left(n_x n_{y} \sqrt{\frac{L_{y y}}{\mu_{y}}}\right)$ \\ [4ex] 
 & \specialcell{Large \\Scale} & $\widetilde{O}\left(\sqrt{\frac{L_{xx}}{\mu_x} +\frac{2L^2_{xy}}{\mu_x\mu_y  }}\right)$ & $\widetilde{O}\left(n_y\sqrt{\frac{L_{xx} L_{yy}}{\mu_x\mu_y} +\frac{2L^2_{xy}}{\mu_x\mu_y} }  \right)$ \\[4ex] 
\hline \hline
\end{tabular}
\end{table}

\section{Solving Min-Min Problems}
Consider the problem
\begin{equation}
    \label{problem:min-min}
    \min_{x\in \mathcal{X}}\min_{y \in \mathbb{R}^{n_y}} f(x, y),
\end{equation}
where $\mathcal{X} \subseteq \mathbb{R}^{n_x}$ is a closed convex set, $f(x,y)$ is a convex function equipped with a mixed oracle, i.e. we have access to a first-order oracle for the outer problem (minimization w.r.t. $x$) and a zeroth-order oracle for the inner problem (minimization w.r.t. $y$). In the sections below we will describe the two approaches to solving such problems together with additional assumptions they require.

The general idea of the approaches is as follows. Let us introduce the function
\begin{equation}
    \label{inner:min-min}
    g(x) = \min_{y\in \mathbb{R}^{n_y}} f(x, y)
\end{equation}
and rewrite the initial problem \eqref{problem:min-min} as
\begin{equation}
    \label{outer:min-min}
    \min_{x\in \mathcal{X}}g(x).
\end{equation}
Using an iterative method for the outer problem \eqref{outer:min-min} requires solving the inner problem \eqref{inner:min-min} numerically on each iteration. An error of the solution of the inner problem results in an inexact oracle for the outer problem.

\subsection{Small dimension of the outer problem}
The approach described in the present subsection requires the following assumptions about the problem \eqref{problem:min-min}:
\begin{enumerate}
    \item $\mathcal{X} \subset \mathbb{R}^{n_x}$ is a compact convex set with nonempty interior;
    \item $n_x$ is relatively small (up to a hundred);
    \item $f(x, y)$ is a continuous convex function which is also $\mu_y$-strongly convex in $y$;
    \item for all $x \in \mathcal{X}$ the function $f(x, \cdot)$ is $L_{yy}$-smooth, i.e.
    \begin{equation*}
        \left\|\nabla_y f(x,y) -\nabla_y f(x,y')\right\|_2  \leq L_{yy} \left\| y - y'\right\|_2\quad \forall y, y' \in \mathbb{R}^{n_y}.
    \end{equation*}
    \item for any $x \in \mathcal{X}$ the minimization problem \eqref{inner:min-min} has solution $y(x)$.
\end{enumerate}
The algorithms used in the proposed approach and related convergence theorems are given in the subsequent paragraphs. Our proposed approach goes as follows:
\begin{approach}\label{appr:first}
    The outer problem \eqref{outer:min-min} is solved via Vaidya's cutting plane method \cite{vaidya1989new,vaidya1996new}. The inner problem \eqref{inner:min-min} is solved via Accelerated Randomized Directional Derivative method for strongly convex functions (ARDDsc) \cite{Dvurechensky_2021}, see Algorithm~\ref{ARDDsc}.\footnote{Here and below instead of ARDDsc we can use Accelerated coordinate descent methods \cite{nesterov2017efficiency,gasnikov2015accelerated} with replacing partial derivatives by finite-differences. In this case we lost opportunity to play on the choice of the norm (that could save $\sqrt{n_y}$-factor in gradient-free oracle complexity estimate \cite{Dvurechensky_2021}), but, we gain a possibility to replace the wort case $L_{yy}$ to the average one (that could be $n_y$-times smaller \cite{nesterov2017efficiency}). At the end this could also save $\sqrt{n_y}$-factor in gradient-free oracle complexity estimate \cite{nesterov2017efficiency}. }
\end{approach}
The complexity of approach \ref{appr:first} is given in the following theorem:
\begin{thm}\label{th:appr_first}
Approach \ref{appr:first} arrives at $\varepsilon$-solution of the problem \eqref{outer:min-min} after\footnote{$\widetilde{O} (\cdot)= O (\cdot)$ up to a small power of logarithmic factor} $\widetilde{O} (n_x)$ calls to the first-order oracle and $\widetilde{O}\left(n_x n_{y} \sqrt{\frac{L_{y y}}{\mu_{y}}}\right)$ calls to the zeroth-order oracle.
\end{thm}
\begin{remark}
As far as the arithmetic complexity of the iteration is concerned, Vaidya's cutting plane method involves inversions of $n_x \times n_x$ matrices, hence the assumption that $n_x$ is relatively small.
\end{remark}
The complexity bounds from theorem \ref{th:appr_first} are derived in the paragraph \textit{Analysis of the approach} which follows the description of algorithms.

\subsubsection{Vaidya's cutting plane method}
Vaidya proposed a cutting plane method \cite{vaidya1989new,vaidya1996new} for solving problems of the form
\begin{equation}\label{problem_vaidya}
    \min_{x \in \mathcal{X}} g(x),
\end{equation}
where $\mathcal{X} \subseteq \mathbb{R}^n$ is a compact convex set with non-empty interior, and $g: \mathcal{X} \to \mathbb{R}$ is a continuous convex function.

We will now introduce the notation and describe the algorithm. Let $P = \{x \in \mathbb{R}^n: \, Ax\geq b\}$ be the bounded full-dimensional polytope, where $A \in \mathbb{R}^{m\times n}$ and $b \in \mathbb{R}^m$. The logarithmic barrier for $P$ is defined as
$$
L(x) := -\sum_{i=1}^{m} \ln \left(a_{i}^{\top} x-b_{i}\right),
$$
where $a_{i}^{\top}$ is the $i^{th}$ row of $A$. The Hessian of $L(x)$ is given by
\begin{equation}\label{hessian}
H(x) =\sum_{i=1}^{m} \frac{a_{i} a_{i}^{\top}}{\left(a_{i}^{\top} x-b_{i}\right)^{2}}
\end{equation}
and is positive definite for all $x$ in the interior of $P$. The \textit{volumetric barrier} for $P$ is defined as
$$
F(x) = \frac{1}{2} \ln \left(\operatorname{det} H(x)\right),
$$
where $\operatorname{det}H(x)$ denotes the determinant of $H(x)$. The point $\bm\omega$ that minimizes $F(x)$ over $P$ will be called the \textit{volumetric center} of $P$.
Let $\sigma_{i}(x)$ be defined as
\begin{equation}\label{volumetric_barrier}
    \sigma_{i}(x)=\frac{a_{i}^{\top} \left(H(x)\right)^{-1} a_{i}}{\left(a_{i}^{\top} x-b_{i}\right)^{2}}, \quad 1 \leq i \leq m.
\end{equation}

Now, let $\mathcal{R}$ be a radius of some Euclidean ball $\mathcal{B}_{\mathcal{R}}$ that contains $\mathcal{X}$. W.l.o.g we will assume that $\mathcal{B}_{\mathcal{R}}$ is centered at the origin. The parameters of the method $\eta>0$ and $\gamma>0$ are small constants such that $\eta \leqslant 10^{-4},$ and $\gamma \leqslant 10^{-3} \eta$. The algorithm starts out with the simplex
\begin{equation}
    P_0=\left\{x \in \mathbb{R}^n: x_{j} \geqslant-\mathcal{R}, j=\overline{1,n},\ \sum_{j=1}^{n} x_{j} \leqslant n \mathcal{R}\right\} \supseteq \mathcal{B}_{\mathcal{R}} \supseteq \mathcal{X}.
\end{equation}
and produces a sequence of pairs $\left(A_k, b_k\right) \in \mathbb{R}^{m_k\times n}\times \mathbb{R}^{m_k}$, such
that the corresponding polytope $P_k = \{x \in \mathbb{R}^n: \, A_k x\geq b_k\}$ always contains a solution of the problem \eqref{problem_vaidya}. At the beginning of each iteration $k$ we have an approximation $z_k$ to the volumetric center of $P_k$ (for more details on computing the approximation see \cite{vaidya1996new,vaidya1989new}). In particular, on the 0-th iteration we can compute the volumetric center explicitly:
\begin{proprep}\label{initial_center}
The volumetric center for $P_0$ is $\bm\omega = \omega \mathbf{1}_n$, where $\omega := \frac{n-1}{n+1} \mathcal{R}$ and $\mathbf{1}_n$ denotes the vector $(1, \dots, 1)^{\top} \in \mathbb{R}^n$.
\end{proprep}
\begin{proof}
$P_0$ can be viewed as $P_0=\left\{x \in \mathbb{R}^n: A_0 x \geq b_0 \right\}$, where $A_0 \in \mathbb{R}^{(n+1) \times n}$ and $b_0 \in \mathbb{R}^{n+1}$ are defined as
\begin{equation}\label{A_and_b}
    b_0 = - \mathcal{R} \left[\begin{array}{cc}
        \mathbf{1}_n \\
        n
    \end{array}\right] \text{ and } A_0 = \left[\begin{array}{cc}
        e_1^{\top} \\
        \vdots \\
        e_n^{\top} \\
        -\mathbf{1}_n^{\top}
    \end{array}\right]
\end{equation}
It is known that the analytic center of simplex (i.e. minimizer of $L(x)$) is also its volumetric center \cite{renegar1988polynomial}. Note that $L(x)$ is a convex function.
\begin{gather*}
    \nabla L(x)=-\sum_{i=1}^{n+1} \frac{a_i}{a_{i}^{\top} x-b_{i}}=-\sum_{i=1}^{n} \frac{e_i}{x_i + \mathcal{R}} + \frac{\mathbf{1}_n}{n\mathcal{R}-\sum_{j=1}^{n} x_j}\\
    \nabla L(\bar{x}) = 0 \iff \frac{1}{n\mathcal{R}-\sum_{j=1}^{n} \bar{x}_{j}} = \frac{1}{\bar{x}_{i} + \mathcal{R}}, i=\overline{1,n}\\
    \bar{x}_{i} = \frac{n-1}{n+1}\mathcal{R},\ i=\overline{1,n}, \text{ and } \bar{x} \in P_0 \Rightarrow \bm\omega = \bar{x}.
\end{gather*}
\end{proof}

For $k\geq 0$, the next polytope  $\left(A_{k+1}, b_{k+1}\right)$ is defined by either adding or removing a constraint to the current
polytope, depending on the values $\left\{\sigma_i(z_k)\right\}_{i=1}^m$ associated to $P_k$:
\begin{enumerate}
    \item If for some $i \in \{1,\ldots, m\}$ one has $\sigma_i(z_k) = \min\limits_{1\leq j\leq m}\sigma_j(z_k) < \gamma$, then $\left(A_{k+1}, b_{k+1}\right)$ is defined by removing the $i^{th}$ row from $\left(A_k, b_k\right)$.
    \item Otherwise, i.e. if $\min\limits_{1\leq j\leq m}\sigma_j(z_k) \geq \gamma$, the oracle is called with the current point $z_k$ as input.
    If $z_k \in \mathcal{X}$, it returns a vector $c_k$, such that $-c_k \in \partial g(z_k)$, i.e. $-c_k$ is a subgradient of $g$ at $z_k$. Otherwise, it returns a vector $c_k$ such that $c_k^{\top} x\geq c_k^{\top} z_k\ \forall x \in \mathcal{X}$. We choose $\beta_k \in \mathbb{R}$ such that $c_k^{\top}z_k \geq \beta_k$ and
    $$
    \frac{c_k^{\top} \left(H(z_k)\right)^{-1} c_k}{\left(c_k^{\top} z_k -\beta_{k}\right)^{2}}=\frac{1}{2} \sqrt{\eta \gamma}.
    $$
    Then we define $\left(A_{k+1}, b_{k+1}\right)$ by adding the row given by $\left(c_k^{\top}, \beta_k\right)$ to $\left(A_k, b_k\right)$.
\end{enumerate}
After $N$ iterations, the method returns a point $x^N := \arg\min\limits_{1 \leq k \leq N} g(z_k)$.

Now, let us introduce the concept of inexact subgradient.
\begin{definition}
    The vector $c \in \mathbb{R}^n$ is called a $\delta$-subgradient of a convex function $g$ at $z \in \operatorname{dom}f$ (we denote $c \in \partial_{\delta}g(z)$), if $$g(x) \geq g(z) + c^{\top} (x - z) -\delta \quad \forall x \in \operatorname{dom}f.$$
\end{definition}
In fact, a $\delta$-subgradient can be used in Vaidya's method instead of the exact subgradient. In this case, we will call the algorithm Vaidya's method with $\delta$-subgradient. We will now present a theorem that justifies this claim.

\begin{thmrep}\label{thm:vaidya}
Let $\mathcal{B}_{\rho}$ and $\mathcal{B}_{\mathcal{R}}$ be some Euclidean balls of radii $\rho$ and $\mathcal{R}$, respectively, such that $\mathcal{B}_{\rho} \subseteq \mathcal{X} \subseteq \mathcal{B}_{\mathcal{R}}$, and let a number $B>0$ be such that $|g(x) - g(x')| \leq B\ \forall x, x' \in \mathcal{X}$. After $N \geq \frac{2n}{\gamma} \ln \left( \frac{n^{1.5} \mathcal{R}}{\gamma \rho} \right) + \frac{1}{\gamma} \ln \pi$ iterations Vaidya's method with $\delta$-subgradient for the problem \eqref{problem_vaidya} returns a point $x^N$ such that
\begin{equation}
    g(x^N) - g(x_*) \leq \frac{B n^{1.5} \mathcal{R}}{\gamma \rho} \exp \left( \frac{\ln \pi -\gamma N}{2n} \right) + \delta,
\end{equation}
where $\gamma>0$ is the parameter of the algorithm and $x_*$ is a solution of the problem \eqref{problem_vaidya}.
\end{thmrep}
\begin{proof}
The proof goes as follows. First, we consider the set $\mathcal{X}^{\varepsilon} := \{(1 - \varepsilon)x_* + \varepsilon x,\ x \in \mathcal{X} \}$ for some $\varepsilon \in [0,1]$. Then, we calculate the number of iterations $N$ sufficient for the volume of the polytope $P_N$ to go below the volume of $\mathcal{X}^{\varepsilon}$. Finally, we show that after $N$ iterations the algorithm arrives at a $(B \varepsilon + \delta)$-suboptimal point.

Fix $\varepsilon \in [0,1]$ and consider the set $\mathcal{X}^{\varepsilon} := \{(1 - \varepsilon)x_* + \varepsilon x, x \in \mathcal{X} \}$. Note that $\mathcal{X}^{\varepsilon} \subseteq P_0$ and
\begin{equation}\label{volume_bound}
    \operatorname{vol} \mathcal{X}^{\varepsilon} = \varepsilon^n \operatorname{vol} \mathcal{X} \geq \varepsilon^n \operatorname{vol} \mathcal{B}_{\rho}.
\end{equation}
The lower bound on the volume of $n$-ball \cite{chen2014inequalities} is
\begin{equation}\label{ball_bound}
    \operatorname{vol} \mathcal{B}_{\rho} \geq \frac{\rho^n}{\sqrt{\pi \left(n+\frac{e}{2}-1 \right)}}\left(\frac{2 \pi e}{n}\right)^{n / 2}.
\end{equation}
We will now derive an upper bound on the volume of polytope $P_N$ at the beginning of the $N$-th iteration.
According to the Vaidya's paper \cite{vaidya1989new}, the following inequality holds:
\begin{equation}\label{ln_vol_bound}
    \ln \left( \operatorname{vol} P_N \right) \leq n \ln \left( \frac{2n}{\gamma} \right) -\rho^{0}-\frac{\gamma}{2} N,
\end{equation}
where $\gamma$ is the parameter of the algorithm and $\rho^{0}$ denotes the minimum of the volumetric barrier at the beginning of the 0-th iteration. To calculate it, we will use proposition \ref{initial_center}, which states that $\bm\omega=\omega \mathbf{1}_{n}$, where $\omega := \frac{n-1}{n+1} \mathcal{R}$.
Observe that
\begin{equation}\label{equalities}
    \mathcal{R} - \omega = \frac{2}{n+1}\mathcal{R},\ \omega + \mathcal{R} = \frac{2n}{n+1}\mathcal{R} = n(\mathcal{R}-\omega).
\end{equation}
Let us calculate $H(\bm\omega)$ and $F(\bm\omega)$:
\begin{gather*}
    H(\bm\omega) \stackrel{\eqref{hessian}}{=} \frac{1}{\left(\omega+\mathcal{R}\right)^{2}} I_n + \frac{1}{\left(n\mathcal{R}-n\omega\right)^{2}} \mathbf{1}_n \mathbf{1}_n^{\top} \stackrel{\eqref{equalities}}{=} \frac{1}{\left(\omega+\mathcal{R}\right)^{2}} \left( I_n + \mathbf{1}_n \mathbf{1}_n^{\top} \right)\\
    \operatorname{det} H(\bm\omega) = \frac{1}{\left(\omega+\mathcal{R}\right)^{2n}} \left( 1 + n \right) = \frac{(n+1)^{2n+1}}{(2n\mathcal{R})^{2n}}.
\end{gather*}
\begin{equation}
    \label{rho_0}
    \rho_0 = F(\bm\omega) = \left( n + \frac{1}{2} \right) \ln (n+1) - n \ln (2n\mathcal{R})
\end{equation}
Using the formulas \eqref{ln_vol_bound} and \eqref{rho_0}, we obtain
\begin{equation}\label{polytope_bound}
    \ln \left( \operatorname{vol} P_N \right) \leq n \ln \left( \frac{4n^2 \mathcal{R}}{\gamma (n+1)} \right) - \frac{1}{2} \ln (n+1) -\frac{\gamma}{2} N.
\end{equation}
Putting together \eqref{volume_bound}, \eqref{ball_bound} and \eqref{polytope_bound}, we will determine a sufficient number of iterations to assure $\operatorname{vol} P_N < \operatorname{vol} \mathcal{X}^{\varepsilon}$:
\begin{gather}
    n \ln \left( \frac{4n^2 \mathcal{R}}{\gamma (n+1)} \right) - \frac{1}{2} \ln (n+1) -\frac{\gamma}{2} N < n \ln \left( \sqrt{\frac{2 \pi e}{n}} \varepsilon \rho \right) - \frac{1}{2} \ln \left( \pi \left(n+\frac{e}{2}-1 \right) \right) \\
    \gamma N > 2n \ln \left( \frac{4n^{2.5} \mathcal{R}}{\gamma (n+1) \sqrt{2 \pi e} \varepsilon \rho} \right) + \ln \left( \frac{\pi \left(n+\frac{e}{2}-1 \right)}{n+1} \right).
\end{gather}
It is sufficient to do
\begin{equation}\label{iterations}
    N = \frac{2n}{\gamma} \ln \left( \frac{n^{1.5} \mathcal{R}}{\gamma \varepsilon \rho} \right) + \frac{1}{\gamma} \ln \pi
\end{equation}
iterations. Since $\operatorname{vol} P_N < \operatorname{vol} \mathcal{X}^{\varepsilon}$, there will be an iteration $j \in \{0, ..., N-1 \}$ and a point $x_{\varepsilon} \in \mathcal{X}^{\varepsilon}$ such that $x_{\varepsilon} \in P_j$ and $x_{\varepsilon} \notin P_{j+1}$. Thus, on the $j$-th iteration we added a new constraint to the polytope.
There exists a point $x \in \mathcal{X}$ such that $x_{\varepsilon} = (1-\varepsilon)x_*+\varepsilon x$, therefore it follows from the convexity of $g$ that
\begin{multline}\label{suboptimal}
    g(x_{\varepsilon}) \leqslant (1-\varepsilon) g(x_*)+\varepsilon g(x) \leqslant (1-\varepsilon) g(x_*)+\varepsilon \left((g(x_*) + B \right) = \\= g(x_*) + B\varepsilon.
\end{multline}
We will now show that the following inequality holds for the query point $z_j$ at the $j$-th iteration:
\begin{equation}\label{delta_gap}
    g(z_j) < g(x_{\varepsilon}) + \delta.
\end{equation}
First, note that $z_j$ belongs to $\mathcal{X}$, otherwise we would cut out a part of $P_j$ that doesn't intersect with $\mathcal{X}$, which is a contradiction since $x_{\varepsilon} \in P_j \setminus P_{j+1}$ and $x_{\varepsilon} \in  \mathcal{X}$. Recall that the constraint $c_j^{\top} x \geq \beta_j$ was added to the polytope on the $j$-th iteration, where $-c_j \in \partial_{\delta}g(z_j)$ and $\beta_j \leq c_j^{\top} z_j$. According to the definition of $\delta$-subgradient,
$$
g(x) \geq g(z_j) - c_j^{\top} (x - z_j) - \delta\quad \forall x \in \mathcal{X}.
$$
Therefore,
\begin{equation*}
    (P_j \setminus P_{j+1}) \cap \mathcal{X} \subseteq \{ x \in \mathcal{X} : c_j^{\top} (x - z_j) < 0 \} \subseteq \{ x \in \mathcal{X} : g(x) > g(z_j) - \delta \},
\end{equation*}
which results in \eqref{delta_gap} since $x_{\varepsilon} \in (P_j \setminus P_{j+1}) \cap \mathcal{X}$. Combining \eqref{iterations}, \eqref{suboptimal} and \eqref{delta_gap}, we conclude the proof.
\end{proof}

\subsubsection{Accelerated Randomized Directional Derivative method}
We refer to the work \cite{Dvurechensky_2021}. For convenience, we present algorithms from this paper, taking into account that the problem will be a classical optimization problem:
\begin{equation}
    \label{opt_problem}
    \min_{x\in \R^{n_x}} f(x).  
\end{equation}
Gradient approximation:
\begin{equation}
    \label{GA}
   \text{grad}_f(x, \tau, \ee) = \frac{n_x}{\tau} \left(f(x+\tau \ee) - f(x)\right)\ee,
\end{equation}
where $\ee \in \mathcal{RS}^{n_x} _{2}(1)$, i.e. be a random vector uniformly distributed on the surface of the unit Euclidean sphere in $\R^{n_x}$.
\begin{definition}
\label{prox func}
Function $d(x): \mathbb{R}^{n_x} \to \mathbb{R}$ is called prox-function if $d(x)$ is $1$-strongly convex w.r.t. $\| \cdot \|_p$-norm and differentiable on $\mathbb{R}^{n_x}$ function. 
\end{definition}
It is worth noting that in the case of $p = 2$, the prox-function $d(x)$ looks like this 
\begin{equation*}
    d(x) = \frac{1}{2}\|x\|^2_2
\end{equation*}
\begin{definition}
\label{Bregman}
Let $d(x): \mathbb{R}^{n_x} \to \mathbb{R}$ is prox-function. For  any  two  points  $x,x' \in \mathbb{R}^{n_x}$ we define Bregman divergence $V_x(x')$ associated with $d(x)$ as follows: 
\begin{equation*}
    V_x(x') = d(x') - d(x) - \langle \nabla d(x), x' - x \rangle.
\end{equation*}
\end{definition}
It is worth noting that in the case of $p = 2$, the Bregman divergence $V_x(x')$ looks like this 
\begin{equation*}
    V_x(x') = \frac{1}{2}\|x' - x\|^2_2
\end{equation*}
Let $x^*$ be fixed point and $x$ be random vector such that $\mathbb{E}_x\|x - x^*\|^2_p \leq R^2_p$, then 
\begin{equation}
    \label{Omega}
    \mathbb{E}_xd\left(\frac{x - x^*}{R_p}\right) \leq \frac{\Omega_p}{2},
\end{equation}
where $\mathbb{E}_x$ denotes the expectation with respect to random vector $x$ and $\Omega_p$ is defined as follows. We note that for $p = 2$ $\Omega_p = 1$

\begin{algorithm}[H]
\caption{Accelerated Randomized Directional Derivative (ARDD) method \cite{Dvurechensky_2021} }
	\label{ARDD}
\begin{algorithmic}
\State 
\noindent {\bf Input:} $x_0$ - starting point, $N$ - number of iterations, $L$ - smoothness parameter, $\tau$.
\State $y_0 := x_0, ~ w_0 = x_0 $
\For {$k=0,1, 2, \ldots, N - 1$ }
    \State Sample $\ee_{k + 1} \in \mathcal{RS}^{n_x} _2(1)$.
    \State Set $$t_k := \frac{2}{k + 2},~ x_{k + 1} := t_k w_k + (1 -t_k)y_k$$
    \State Calculate $\text{grad}_f(x_{k + 1}, \tau, \ee_{k + 1})$ using \eqref{GA} 
    \State Compute $$y_{k + 1} := x_{k+1} - \frac{1}{2L}\text{grad}_f(x_{k + 1}, \tau, \ee_{k + 1})$$
    \State Set $$\alpha_k := \frac{k + 1}{96n^2L}$$
    \State Compute $$w_{k + 1} := \argmin\limits_{z \in \R^{n_x}}\left\{\alpha_{k+1}\langle \text{grad}_f(x_{k + 1}, \tau_k, \ee_{k + 1}), w - w_k\rangle + V_{z_k}(z)\right\} $$
\EndFor
\State 
\noindent {\bf Output:} $y_{N}$ or $\bar y_{N}$.
\end{algorithmic}
\end{algorithm}

\begin{algorithm}[H]
\caption{Accelerated Randomized Directional Derivative method for strongly convex functions (ARDDsc) \cite{Dvurechensky_2021}}
	\label{ARDDsc}
\begin{algorithmic}
\State
\noindent {\bf Input:} $x_0$ - starting point s.t. $\|x_0 - x^*\|^2_{p}$, $N$ -number of iterations, $\mu_p$ -strong convexity parameter.
\State Set 
    \begin{equation*}
        N_0 = \left\lceil\sqrt{\frac{8aL_2\Omega_p}{\mu_p}}\right\rceil
    \end{equation*}
    where $a = 384n^2\rho_n$, $\rho_n = \min\left\{q-1, 16 \text{ln}n - 8\right\} n^{\frac{2}{q} - 1}$.
\For {$k=0,1, 2, \ldots, N - 1$ }
\State Set 
    $R^2_k = R^2_p2^{-k} $ 
\State Set $d_k(x) = R^2_k d\left(\frac{x - u_k}{R_k}\right)$
\State Run ARDD \ref{ARDD} with starting point $u_k$ and prox-function $d_k(x)$ \ref{prox func} for $N_0$ steps.
\State Set $u_{k + 1} = y_{N_0},~ k = k+1$
\EndFor
\State 
\noindent {\bf Output:} $u_{N}$ 
\end{algorithmic}
\end{algorithm}
\begin{thm}[see \cite{Dvurechensky_2021}]
\label{ARDDsc convergence}
Let $p \in [1, 2]$ and $q \in [2, +\infty]$ be defined such that $\frac{1}{p} + \frac{1}{q} = 1$.
Let function $f$ in problem \eqref{opt_problem} be $\mu_p$-strongly convex w.r.t. $\|.\|_p$ and $L_2$-smooth w.r.t. $\|.\|_2$ and ARDDsc method \ref{ARDDsc} be applied to
solve this problem. Then 
\begin{equation*}
    \mathbb{E}f(u_N) - \min_{x \in \mathbb{R}^{n_x}}f(x) \leq \frac{\mu_pR^2_p}{2}2^{-N} 
\end{equation*}
Moreover, the oracle complexity to achieve $\varepsilon$-accuracy of the solution is
\begin{equation*}
    \widetilde{O}\left(n^{\frac{1}{q}+ \frac{1}{2}}\sqrt{\frac{L_2\Omega_p}{\mu_p}}\log_2\frac{\mu_p R^2_p}{\varepsilon}\right)
\end{equation*}
\end{thm}
\subsubsection{Analysis of the approach \ref{appr:first}}
Fix a point $x' \in \mathcal{X}$. The following theorem gives the recipe to obtaining the $\delta$-subgradient $c \in \partial_{\delta} g(x')$ for the outer problem \eqref{outer:min-min}:
\begin{thmrep}\label{thm:delta_min_min}
    Let $\tilde{\delta}>0$ and $\tilde{y} \in \R^{n_y}$ satisfy $f(x', \tilde{y}) - g(x') \leq \tilde{\delta}$, then $\partial_{x} f(x', \tilde{y}) \in \partial_{\delta} g(x')$ with
    \begin{equation}\label{delta_min_min}
        \delta = 2 \left(\sqrt{\tilde{\delta}} + 2\sqrt{D} \right) \sqrt{\frac{L_{yy}\tilde{\delta}}{\mu_{y}}},
    \end{equation}
    where $D:= \displaystyle \max_{x \in \mathcal{X}} \left( f(x, \mathbf{0}) - g(x) \right)<+\infty$.
\end{thmrep}
\begin{proof}
Using Lemmas \ref{bound_dot_prod}, \ref{statement_gas} and triangle inequality, we get $\partial_{x} f(x', \tilde{y}) \in \partial_\delta g(x')$ with
\begin{equation}\label{delta_theorem}
    \delta = \left(\|\tilde{y}\|_2 + \max_{x \in \mathcal{X}}\|y(x)\|_2 \right) \sqrt{2L_{yy} \tilde{\delta}}.
\end{equation}
Let us bound the two terms in parentheses. For any $x \in \mathcal{X}$, it follows from the strong convexity of $f(x, \cdot)$ that
\begin{equation}\label{y_deviation_bound}
    \|y' - y(x)\|_2 \leq \sqrt{\frac{2(f(x, y') - g(x))}{\mu_{y}}}\quad \forall y' \in \R^{n_y}.
\end{equation}
In particular, putting $y'=\mathbf{0}$ leads to
\begin{equation*}
    \|y(x)\|_2 \leq \sqrt{\frac{2(f(x, \mathbf{0}) - g(x))}{\mu_{y}}}.
\end{equation*}
Observe that $\Delta(x) := f(x, \mathbf{0}) - g(x)$ is a continuous function defined on the compact set $\mathcal{X}$. Therefore, it is bounded, i.e. $\Delta(x) \leq D$ for some $D<+\infty$. Thus,
\begin{equation}\label{solution_bound}
    \max_{x \in \mathcal{X}}\|y(x)\|_2 \leq \sqrt{\frac{2D}{\mu_{y}}}.
\end{equation}
Finally, let us bound $\| \tilde{y} \|_2$ by again using triangle inequality, then inequalities \eqref{y_deviation_bound}, \eqref{solution_bound} and $f\left(x^{\prime}, \tilde{y}\right)-g\left(x^{\prime}\right) \leq \tilde{\delta}$:
\begin{equation}\label{y_tilde_bound}
    \| \tilde{y} \|_2 \leq \| \tilde{y} - y(x') \|_2 + \| y(x') \|_2 \leq \| \tilde{y} - y(x') \|_2 + \sqrt{\frac{2D}{\mu_{y}}} \leq \sqrt{\frac{2 \tilde{\delta}}{\mu_{y}}} + \sqrt{\frac{2D}{\mu_{y}}}.
\end{equation}
The statement of the theorem immediately follows from \eqref{delta_theorem}, \eqref{solution_bound} and \eqref{y_tilde_bound}.
\end{proof}
The theorem \ref{thm:delta_min_min} is based on the two following lemmas:
\begin{lmmrep}\label{bound_dot_prod}
    Let $h: \R^{n_y} \to \R$ be an $L$-smooth convex function, and let the point $\tilde{y} \in \R^{n_y}$ satisfy $h(\tilde{y}) - h(y_*) \leq \tilde{\delta}$ for some $\tilde{delta}>0$, where $y_* \in \underset{y \in \R^{n_y}}{\operatorname{Argmin}}\ h(y)$. Then
    \begin{equation*}
        \left\langle \nabla h(\tilde{y}), \tilde{y}-y\right\rangle \leq \left\| \tilde{y} - y \right\|_2 \sqrt{2L \tilde{\delta}}\quad \forall y \in \mathbb{R}^{n_y}.
    \end{equation*}
\end{lmmrep}
\begin{proof}
Due to the Cauchy–Bunyakovsky–Schwarz inequality,
\begin{equation}\label{bunyakovsky}
    \left\langle \nabla h(\tilde{y}), \tilde{y}-y\right\rangle \leq \left\|  \nabla h(\tilde{y}) \right\|_2 \left\| \tilde{y}-y \right\|_2
\end{equation}
Since $h$ is $L$-smooth and convex,
\begin{equation*}
    h(\tilde{y}) \geq h(y_*) + \langle \nabla h(y_*), \tilde{y} - y_*\rangle + \frac{1}{2L} \| \nabla h(y_*) - \nabla h(\tilde{y}) \|_2^2
\end{equation*}
Using $\nabla h(y_*)=0$, we obtain
\begin{equation}\label{grad_norm_bound}
    \| \nabla h(\tilde{y}) \|_2^2 \leq 2L \left( h(\tilde{y}) - h(y_*) \right).
\end{equation}
Combining \eqref{bunyakovsky} and \eqref{grad_norm_bound}, we conclude the proof.
\end{proof}
\begin{lmmrep}[see \cite{gas2015optima}, p.12]\label{statement_gas}
    Let $\tilde{y} \in \R^{n_y}$ satisfy
    \begin{equation}\label{dot_prod_bound}
        \left\langle \nabla_y f(x, \tilde{y}), \tilde{y}-y(x')\right\rangle \leq \delta\quad \forall x' \in \mathcal{X},
    \end{equation}
    then $\partial_{x} f(x, \tilde{y}) \in \partial_\delta g(x)$.
\end{lmmrep}
\begin{proof}
From convexity of $f$ we get for any $x' \in \mathcal{X}$
\begin{equation*}
    f(x', y(x')) \geq f(x, \tilde{y}) + \langle \nabla_x f(x, \tilde{y}), x' - x \rangle + \langle \nabla_y f(x, \tilde{y}), y(x') - \tilde{y} \rangle
\end{equation*}
Using \eqref{dot_prod_bound} and $g(x') = f(x', y(x'))$, we obtain
\begin{equation*}
    g(x') \geq f(x, \tilde{y}) + \langle \nabla_x f(x, \tilde{y}), x' - x \rangle - \delta.
\end{equation*}
Note that $g(x) \leq f(x, \tilde{y})$, therefore,
\begin{equation*}
    g(x') \geq g(x) + \langle \nabla_x f(x, \tilde{y}), x' - x \rangle - \delta\quad \forall x' \in \mathcal{X}.
\end{equation*}
Thus, by definition we have $\partial_{x} f(x, \tilde{y}) \in \partial_\delta g(x)$.
\end{proof}
According to theorem \ref{thm:delta_min_min}, we need to solve the inner problem \eqref{inner:min-min} with sufficient accuracy to obtain the $\delta$-subgradient. In fact, we can simplify the formula \eqref{delta_min_min} to
\begin{equation}\label{delta_simple}
    \delta = 6 \sqrt{\frac{L_{yy}D\tilde{\delta}}{\mu_{y}}}.
\end{equation}
Indeed, since $\tilde{y}$ is an $\tilde{\delta}$-solution of the inner problem \eqref{inner:min-min}, we can assume that $\tilde{\delta} \leq D \equiv \max_{x \in \mathcal{X}}(f(x, \mathbf{0})-g(x))$. If this inequality doesn't hold, we can always take $\mathbf{0}$ as an approximate solution of \eqref{inner:min-min}.

Now, to derive the complexity of approach \ref{appr:first}, we will use the theorem \ref{thm:vaidya} and put
\begin{equation*}
    \frac{B n^{1.5} \mathcal{R}}{\gamma \rho} \exp \left(\frac{\ln \pi-\gamma N}{2 n}\right) = \frac{\varepsilon}{2} \text{ and } \delta = \frac{\varepsilon}{2},
\end{equation*}
i.e. Vaidya's method will perform
\begin{equation*}
    N_x = O \left( n_x \ln \left( \frac{n_x^{1.5} B \mathcal{R}}{\varepsilon \rho} \right) \right),
\end{equation*}
steps (first-order oracle calls), and at each of them ARDDsc will perform
\begin{equation*}
    N_y = \widetilde{O} \left( n_{y} \sqrt{\frac{L_{y y}}{\mu_{y}}} \right)
\end{equation*}
iterations (see theorem \ref{ARDDsc convergence}). Thus, the number of zeroth-order oracle calls is
\begin{equation*}
    N_x \cdot N_y = \widetilde{O}\left(n_x n_{y} \sqrt{\frac{L_{y y}}{\mu_{y}}}\right),
\end{equation*}
which finishes the analysis of approach \ref{appr:first}.



\section{Solving Min-Max Saddle-Point Problems}
Consider the problem
\begin{equation}
    \label{problem}
    \min_{x\in \mathcal{X}}\max_{y \in \mathbb{R}^{n_y}} f(x, y),
\end{equation}
where $\mathcal{X} \subseteq \mathbb{R}^{n_x}$ is a closed convex set, $f(x,y)$ is a convex-concave function (i.e. convex in $x$ and concave in $y$) equipped with a mixed oracle, i.e. we have access to a first-order oracle for the outer problem (minimization w.r.t. $x$) and a zeroth-order oracle for the inner problem (maximization w.r.t. $y$). In the subsections below we describe the two approaches for solving such problems together with additional assumptions they require.

The general idea of the approaches is as follows. Let us introduce the function
\begin{equation}
    \label{inner}
    g(x) = \max_{y\in \mathbb{R}^{n_y}} f(x, y)
\end{equation}
and rewrite the initial problem \eqref{problem} as
\begin{equation}
    \label{outer}
    \min_{x\in \mathcal{X}}g(x).
\end{equation}
Using an iterative method for the outer problem \eqref{outer} requires solving the inner problem \eqref{inner} numerically in each iteration. An error of the solution of the inner problem results in an inexact oracle for the outer problem.



\subsection{Small dimension of the outer problem}
The approach described in the present section requires the following assumptions about the problem \eqref{problem}:
\begin{enumerate}
    \item $\mathcal{X} \subset \mathbb{R}^{n_x}$ is a compact convex set with nonempty interior;
    \item $n_x$ is relatively small (up to a hundred);
    \item $f(x, y)$ is a continuous function which is convex in $x$ and $\mu_y$-strongly concave in $y$;
    \item for all $x \in \mathcal{X}$ the function $f(x, \cdot)$ is $L_{yy}$-smooth, i.e.
    \begin{equation*}
        \left\|\nabla_y f(x,y) -\nabla_y f(x,y')\right\|_2  \leq L_{yy} \left\| y - y'\right\|_2\quad \forall y, y' \in \mathbb{R}^{n_y}.
    \end{equation*}
    \item for any $x \in \mathcal{X}$ the maximization problem \eqref{inner} has solution $y(x)$.
\end{enumerate}
The algorithms used in the approach and related convergence theorems are given in the previous section. The approach goes as follows:
\begin{approach}\label{appr:first_minmax}
    The outer problem \eqref{outer} is solved via Vaidya's cutting plane method. The inner problem \eqref{inner:min-min} is solved via 
    ARDDsc, see Algorithm~\ref{ARDDsc}.
\end{approach}
Complexity of the approach is given in the following theorem:
\begin{thm}\label{th:appr_first_minmax}
Approach \ref{appr:first_minmax} arrives at $\varepsilon$-solution of the problem \eqref{outer} after $\widetilde{O} (n_x)$ calls to the first-order oracle and $\widetilde{O}\left(n_x n_{y} \sqrt{\frac{L_{y y}}{\mu_{y}}}\right)$ calls to the zeroth-order oracle.
\end{thm}
\begin{remark}
As far as the arithmetic complexity of the iteration is concerned, Vaidya's cutting plane method involves inversions of $n_x \times n_x$ matrices, hence the assumption that $n_x$ is relatively small.
\end{remark}
Complexity bounds from theorem \ref{th:appr_first_minmax} are derived in the following paragraph.

\subsubsection{Analysis of the approach \ref{appr:first_minmax}}
Fix a point $x' \in \mathcal{X}$. The following lemma gives the recipe to obtaining the $\delta$-subgradient $c \in \partial_{\delta} g(x')$ for the outer problem \eqref{outer}:
\begin{lmm}[see \cite{polyak1983intro}]\label{lem:delta_min_max}
Let $\tilde{y} \in \mathbb{R}^{n_{y}}$ satisfy $g\left(x^{\prime}\right) - f\left(x^{\prime}, \tilde{y}\right) \leq \delta,$ then $\partial_{x} f\left(x^{\prime}, \tilde{y}\right) \in \partial_{\delta} g\left(x^{\prime}\right)$.
\end{lmm}
According to lemma \ref{lem:delta_min_max}, we need to solve the inner problem \eqref{inner} with accuracy $\delta$ to obtain the $\delta$-subgradient.

Now, to derive the complexity of approach \ref{appr:first_minmax}, we will use the theorem \ref{thm:vaidya} and put
\begin{equation*}
    \frac{B n^{1.5} \mathcal{R}}{\gamma \rho} \exp \left(\frac{\ln \pi-\gamma N}{2 n}\right) = \frac{\varepsilon}{2} \text{ and } \delta = \frac{\varepsilon}{2},
\end{equation*}
i.e. Vaidya's method will perform
\begin{equation*}
    N_x = O \left( n_x \ln \left( \frac{n_x^{1.5} B \mathcal{R}}{\varepsilon \rho} \right) \right),
\end{equation*}
steps (first-order oracle calls), and at each of them ARDDsc will perform
\begin{equation*}
    N_y = \widetilde{O} \left( n_{y} \sqrt{\frac{L_{y y}}{\mu_{y}}} \right)
\end{equation*}
iterations (see theorem \ref{ARDDsc convergence}). Thus, the number of zeroth-order oracle calls is
\begin{equation*}
    N_x \cdot N_y = \widetilde{O}\left(n_x n_{y} \sqrt{\frac{L_{y y}}{\mu_{y}}}\right),
\end{equation*}
which finishes the analysis of approach \ref{appr:first_minmax}.

\subsection{Large dimension of the outer problem}

For a detailed study of the convergence of the methods, we introduce some assumptions about the objective function $f(x, y)$.

\textbf{Assumption 1.} $f(x,y)$ is convex-concave. It means that $f(\cdot, y)$ is convex for all $y$ and  $f(x, \cdot)$ is concave for all $x$.

\textbf{Assumption 1(s).} $f(x,y)$ is strongly-convex-strongly-concave. It means that $f(\cdot, y)$ is $\mu_x$-strongly convex for all $y$ and  $f(x, \cdot)$ is $\mu_y$-strongly concave for all $x$ w.r.t. $\|\cdot\|_2$, i.e. for all $x_1, x_2 \in \mathcal{X}$ and for all $y_1, y_2 \in \mathbb{R}^{n_y}$ we have
\begin{eqnarray}
\label{SC}
f(x_1,y_2) &\geq& f(x_2, y_2) + \langle\nabla_x f(x_2, y_2) , x_1 - x_2\rangle  + \frac{\mu_x}{2}\|x_1 - x_2\|^2_2, \nonumber\\
-f(x_2,y_1) &\geq& -f(x_2, y_2) - \langle\nabla_y f(x_2, y_2) , y_1 - y_2\rangle + \frac{\mu_y}{2}\|y_1 - y_2\|^2_2 .
\end{eqnarray}

\textbf{Assumption 2.} $f(x,y)$ is $(L_{xx}, L_{xy}, L_{yy})$-smooth w.r.t $\|\cdot\|_2$, i.e. for all $x, x' \in \mathcal{X},~ y, y' \in \mathbb{R}^{n_y}$ 
\begin{eqnarray}
\label{smooth}
    \left\|\nabla_x f(x,y) -\nabla_x f(x',y)\right\|_2 &\leq& L_{xx} \left\| x - x'\right\|_{2};\nonumber \\ \left\|\nabla_x f(x,y) -\nabla_x f(x,y')\right\|_{2}  &\leq& L_{xy} \left\| y - y'\right\|_{2} \nonumber\\ \left\|\nabla_y f(x,y) -\nabla_y f(x',y)\right\|_{2}  &\leq& L_{xy} \left\| x - x'\right\|_{2}.\nonumber\\
    \left\|\nabla_y f(x,y) -\nabla_y f(x,y')\right\|_{2}  &\leq& L_{yy} \left\| y - y'\right\|_{2}.
\end{eqnarray}

As mentioned above, we have access to a first-order oracle $\nabla_x f(x, y)$ for the outer problem (minimization problem with variables $x$) and a zeroth-order oracle $f(x, y)$ for the inner problem (maximization problem with variables $y$). Since we do not have access to the values of the gradient $\nabla_y f(x, y)$, it is logical to approximate it using finite differences using the value of the function $f(x, y)$ at two close points as follows
\begin{equation}
    \label{approximation}
    \text{grad}_f(x, y, \tau, \ee) = -\frac{n}{\tau}\left(f(x, y + \tau\ee) - f(x, y)\right) \ee,
\end{equation}
where $\ee \in \mathcal{RS}^{n_y} _{2}(1)$, i.e. is a random vector uniformly distributed on the surface of the unit Euclidean sphere in $\R^{n_y}$. So we get a mixed oracle
\begin{equation}
\label{mixed oracle}
   G(x, y, \tau, \ee) = \begin{pmatrix} \nabla_x f(x,y)\\
                        \text{grad}_f(x, y, \tau, \ee)
\end{pmatrix}.
\end{equation}

Using the mixed oracle \eqref{mixed oracle}, we provide our approach for solving the initial saddle-point problem \eqref{problem}. First, we can use the following trick with the help of Sion's theorem:
\begin{equation*}
    \min_{x\in \mathcal{X}}\max_{y\in \mathbb{R}^{n_y}} f(x, y)  = \max_{y\in \mathbb{R}^{n_y}} \min_{x\in \mathcal{X}}  f(x, y) =  \max_{y\in \mathbb{R}^{n_y}} h(y), \text{ where } h(y) = \min_{x\in \mathcal{X}} f(x, y).
\end{equation*}
For the new problem, we apply the Catalyst algorithm \cite{catalyst} to the outer maximization problem:
\begin{equation}
\label{outer catalyst}
     \max_{y\in \mathbb{R}^{n_y}}\left\{ h(y) = \min_{x\in \mathcal{X}} f(x, y)\right\}.
\end{equation}
\begin{algorithm}[H]
\caption{Catalyst \cite{catalyst}}
\label{Catalyst}
\begin{algorithmic}[1]
\State 
\noindent {\bf Input:} starting point $x_0$, parameters $H_1$ and $\alpha_0$, accuracy of solution to subproblem  $\tilde{\varepsilon}$, optimization method $\mathcal{M}$.
\State Initialize $$q = \frac{\mu_y}{(\mu_y + H_1)}$$
\While {the desired stopping criterion is not satisfied}
    \State Find an approximate solution of the following problem using  $\mathcal{M}$: 
    \begin{eqnarray}
    \label{subpr}
         y_k &\approx& \argmax_{y \in \mathbb{R}^{n_y}}\left\{ \varphi_k(y) = h(y) - \frac{H_1}{2}\|y - z_{k - 1}\|^2_2\right\},\\&&\text{ such that } \varphi^*_k - \varphi_k(y_k) \leq \tilde{\varepsilon} \nonumber
    \end{eqnarray}
    \State Compute $\alpha_k \in (0, 1)$ from equation $$\alpha^2_k = (1- \alpha_i)\alpha^2_{k - 1} + q\alpha_k $$
    \State Compute $$z_k = y_k + \beta_k\left(y_k - y_{k - 1}\right), \text{ where } \beta_k = \frac{\alpha_{k - 1}(1 - \alpha_{k - 1})}{\alpha^2_{k - 1}+ \alpha_k}$$
        
\EndWhile 
\State
\noindent {\bf Output:} $y_{final}$.
\end{algorithmic}
\end{algorithm}
Now the question arises how to solve the auxiliary problem \eqref{subpr}. This subproblem is  equivalent to solving the following problem:
\begin{equation*}
    \max_{y\in \mathbb{R}^{n_y}}\left\{ h(y) - \frac{H_1}{2}\|y - z_{k - 1}\|^2_2\right\} = \max_{y\in \mathbb{R}^{n_y}}\left\{ \min_{x \in \mathcal{X}}\left\{ f(x, y) \right\} - \frac{H_1}{2}\|y - z_{k - 1}\|^2_2 \right\},
\end{equation*}
for which we again use the Sion's theorem and rewrite equivalently the problem as:
\begin{equation*}
    \max_{y \in \mathbb{R}^{n_y}} \min_{x\in \mathcal{X}}\left\{ f(x, y) - \frac{H_1}{2}\|y - z_{k - 1}\|^2_2 \right\} = \min_{x \in \mathcal{X}}\max_{y \in \mathbb{R}^{n_y}}\left\{ f(x, y)  - \frac{H_1}{2}\|y - z_{k - 1}\|^2_2 \right\}.
\end{equation*}
For convenience, we  denote
\begin{equation}
\label{catalyst function}
    \psi(x, y) = f(x, y) - \frac{H_1}{2}\|y - z_k\|^2_2.
\end{equation}
Thus, to solve the auxiliary problem \eqref{subpr}, we first solve the following saddle-point problem
\begin{equation}
\label{catalyst problem}
    \min_{x \in \mathcal{X}}\max_{y \in \mathbb{R}^{n_y}} \psi(x, y) = \min_{x \in \mathcal{X}}\max_{y \in \mathbb{R}^{n_y}}\left\{ f(x, y)  - \frac{H_1}{2}\|y - z_{k - 1}\|^2_2 \right\}.
\end{equation}
This saddle-point problem \eqref{catalyst problem} can be considered as an optimization problem for a certain function. Indeed, let  us introduce a function
\begin{equation}
    \label{catalyst inner}
    \xi(x) = \max_{y\in \mathbb{R}^{n_y}} \psi(x, y),
\end{equation}
and rewrite the initial problem \eqref{catalyst problem} as follows:
\begin{equation}
    \label{catalyst outer}
    \min_{x\in \mathcal{X}}\xi(x).
\end{equation}
To solve problem \eqref{catalyst problem}, we solve the outer minimization problem with respect to the variable $x$ by the fast adaptive gradient method with inexact oracle. In each iteration of this method, to find the inexact first-order oracle for the outer problem, we solve the inner problem \eqref{catalyst inner}. Since for this inner problem we have access only to the zeroth-order oracle, we use accelerated gradient-free method ARDDsc \cite{Dvurechensky_2021}.
Our approach is summarized as follows.
\begin{approach}
\label{large scale minimax approach}
    The outer problem \eqref{outer} is solved via Catalyst Algorithm \ref{Catalyst}. The subproblem \eqref{subpr}    is solved as saddle-point problem \eqref{catalyst problem}. The outer problem \eqref{catalyst outer} is solved via Fast Adaptive Gradient Method (Algorithm \ref{Fast adaptive gradient method}). At each iteration of Algorithm \ref{Fast adaptive gradient method} the inner problem \eqref{catalyst inner} is solved via ARDDsc, see Algorithm~\ref{ARDDsc},
    for case $p = 2$ (that is, prox-function $d(x) = \frac{1}{2}\|x\|^2_2$, see Definition \ref{prox func}, Bregman divergence $V_x(x')= \frac{1}{2}\|x' - x\|^2_2$, see Definition \ref{Bregman}, and $\Omega_p = \Omega_2 = 2$).
\end{approach}

\subsubsection{Analysis of fast adaptive gradient method with $(\delta, \sigma, L, \mu)$-oracle (Algorithm \ref{Fast adaptive gradient method})}
For our analysis, due to the fact that Algorithm \ref{ARDDsc} is randomized, we need not just a fast gradient method for solving the outer problem \eqref{catalyst outer}, but a fast adaptive gradient method for $(\delta, \sigma, L, \mu)$-oracle. This is the extension of the fast adaptive gradient method from \cite{stonyakin2020inexact}. To understand this problem in depth and in detail, we need to carefully consider the concept of $(\delta, \sigma, L, \mu)$- oracle and perform a deep and thorough convergence analysis of the fast gradient method using such a seemingly unusual oracle. To that end, we consider the following general minimization problem:
\begin{equation}
\label{fgd_pr}
    \min_{x\in\mathcal{X}}f(x).
\end{equation}
\begin{definition}[$(\delta, \sigma, L, \mu)$-oracle]
Let function f be convex on convex set $\mathcal{X}$. We say that it is equipped with a first-order $(\delta, \sigma, L, \mu)$-oracle if, for any $x' \in \mathcal{X}$, we can compute a pair $\left(f_{\delta, L, \mu}(x'), g_{\delta, L, \mu}(x')\right) \in \mathbb{R}\times\mathbb{R}^{n_x}$ such that with probability at least $1-\sigma$
\begin{equation}
\label{delta sigma oracle}
    \frac{\mu}{2}\|x' - x\|^2 \leq f(x) - \left(f_{\delta, L, \mu}(x') +\left\langle g_{\delta, L, \mu}(x'), x- x' \right\rangle\right) \leq \frac{L}{2}\|x - x'\|^2 + \delta.
\end{equation}
\end{definition}

\begin{definition}[$(\varepsilon, \sigma)$-solution]
\label{es solution}
Let $\varepsilon > 0$ be the target accuracy of the solution and $\sigma \in (0, 1)$ be the target confidence level. We say that a random point $\hat{x} \in \mathcal{X}$ is $(\varepsilon, \sigma)$-solution to problem \eqref{fgd_pr} if
\begin{equation}
\label{delta sigma solution}
    \mathbb{P}\left\{f(\hat{x}) - \min_{x\in\mathcal{X}}f(x) \leq \varepsilon\right\} \geq 1-\sigma.
\end{equation}
If $\sigma = 0$, we say that $\hat{x} \in \mathcal{X}$ is an $\varepsilon$-solution to problem \eqref{fgd_pr}.
\end{definition}

\begin{algorithm}[H]
\caption{Fast adaptive gradient method with $(\delta, \sigma, L, \mu)$-oracle \cite{stonyakin2020inexact}}
\label{Fast adaptive gradient method}
\begin{algorithmic}[1]
\State 
\noindent {\bf Input:} starting point $x_0$, $L_0 > 0$, $\mu \geq 0$, sequence $\{\delta\}_{k\geq 0}$.
\State $y_0:= x_0,~ u_0:= x_0,~ \alpha_0 := 0,~ A_0 := \alpha_0$
\For {$k \geq 0$}
    \State Find the smallest integer  $i_k \geq 0$ such that
    \begin{equation*}
        f_{\delta_k, L, \mu}(x_{k + 1}) \leq  f_{\delta_k, L, \mu}(y_{k + 1}) + \left\langle g_{\delta,L, \mu}(y_{k + 1}), x_{k + 1} - y_{k + 1}\right\rangle + \frac{L_{k + 1}}{2}\|x_{k + 1} - y_{k + 1}\|^2_{2} + \delta_k,
    \end{equation*}
    where $L_{k + 1} = 2^{i_k - 1}L_k$.
    \State  Compute $\alpha_{k + 1}$ such that $\alpha_{k + 1}$ is the largest root of
    \begin{equation*}
        A_{k + 1}(1 + A_k\mu) = L_{k + 1} \alpha^2_{k + 1}, \text{ where } A_{k+1}:= A_k + \alpha_{k + 1} 
    \end{equation*}
    \State $y_{k+1} = \frac{\alpha_{k + 1}u_k + A_k x_k}{A_{k +1}}$
    \begin{equation*}
        \phi_{k + 1}(x) = \alpha_{k + 1}\left\langle g_{\delta,L, \mu}(y_{k + 1}), x - y_{k + 1}\right\rangle + \frac{(1 +  A_k\mu)}{2}\|x - u_k\|^2_2 + \frac{\alpha_{k+1}\mu}{2}\|x  - y_{k + 1}\|^2_2
    \end{equation*}
    \State $u_{k+1}:= \argmin\limits_{x \in \mathcal{X}} \phi_{k + 1}(x)$
    \State $x_{k+1} = \frac{\alpha_{k + 1}u_{k+1} + A_k x_k}{A_{k +1}}$
\EndFor
\State
\noindent {\bf Output:} $x_{k + 1}$.
\end{algorithmic}
\end{algorithm}

Note that the problem $\argmin\limits_{x \in \mathcal{X}} \phi_{k + 1}(x)$ is solved exactly in each iteration.

\begin{thmrep}
\label{fgm theorem}
Let function $f$ be convex on convex set $\mathcal{X}$ and be  equipped with a first-order $(\delta, \sigma, L, \mu)$-oracle.
Then, after $N$ iterations of Algorithm \ref{Fast adaptive gradient method} applied to problem \eqref{fgd_pr}, we have that with probability at least $(1 - N\sigma)$:
\begin{eqnarray}
\label{fgm convergence}
    f(x_N) - f(x^*)&\leq&2L\exp\left(-\frac{N - 1}{2}\sqrt{\frac{\mu}{L}}\right) R^2_2 + \frac{2\sum^{N-1}_{k = 0}A_{k + 1} \delta_k}{A_N},
\end{eqnarray}
where $R_2$ is such that $\frac{1}{2}\|x_0 - x^*\|^2_2 \leq R^2_2$ and $x_0$ is the starting point.
\end{thmrep}
\begin{proof}
We almost completely repeat the proof of result from \cite{stonyakin2020inexact}.
\begin{lemma}[see \cite{stonyakin2020inexact}]
Let $\phi(x)$ be convex function on the convex set $\mathcal{X}$ and 
\begin{equation*}
    y = \argmin_{x \in \mathcal{X}}\left\{\phi(x) + \frac{\beta}{2}\|x -z\|^2_2 + \frac{\gamma}{2}\|x - u\|^2_2\right\},
\end{equation*}
where $\beta \geq 0$, $\gamma \geq 0$. Then 
\begin{equation}
\label{lemma stonyakin}
    \psi(x) + \frac{\beta}{2}\|x -z\|^2_2 + \frac{\gamma}{2}\|x - u\|^2_2 \geq  \psi(y) + \frac{\beta}{2}\|y -z\|^2_2 + \frac{\gamma}{2}\|y - u\|^2_2 + \frac{\gamma + \beta}{2}\|y - x\|^2_2
\end{equation}
\end{lemma}

\begin{lemma}[see \cite{stonyakin2020inexact}]
For all $x \in \mathcal{X}$ we have  with probability with $1 -\sigma$
\begin{eqnarray}
\label{convergence lemma}
    A_{k+1}f(x_{k +1}) - A_{k}f(x_{k}) &+& \frac{(1 +  A_{k+1}\mu)}{2}\|x - u_{k+1}\|^2_2 - \frac{(1 +  A_k\mu)}{2}\|x - u_k\|^2_2 \nonumber\\ &\leq& \alpha_{k+1} f(x) + 2\delta_{k+1}A_{k+1}
\end{eqnarray}
\end{lemma}

\textbf{Proof} Using definition of  $(\delta, \sigma, L, \mu)$-oracle \eqref{delta sigma oracle}, we have with probability with $1 -\sigma$ inequalities $f(x) - \delta \leq f_{\delta, L, \mu}(x) \leq f(x)$. Then with probability at least $1-\sigma$:
\begin{equation*}
    f(x_{k + 1}) \leq  f_{\delta_k, L, \mu}(y_{k + 1}) + \left\langle g_{\delta,L, \mu}(y_{k + 1}), x_{k + 1} - y_{k + 1}\right\rangle + \frac{L_{k + 1}}{2}\|x_{k + 1} - y_{k + 1}\|^2_{2} + 2\delta_k,
\end{equation*}
Using definitions of sequences $x_{k +1}$ and $y_{k +1}$, we have with probability at least $1-\sigma$:
\begin{eqnarray*}
    f(x_{k + 1}) &\leq&  f_{\delta_k, L, \mu}(y_{k + 1}) + \left\langle g_{\delta,L, \mu}(y_{k + 1}), x_{k + 1} - y_{k + 1}\right\rangle + \frac{L_{k + 1}}{2}\|x' - x\|^2_{2} + 2\delta_k\\ &=& f_{\delta_k, L, \mu}(y_{k + 1}) + \left\langle g_{\delta,L, \mu}(y_{k + 1}), \frac{\alpha_{k + 1}u_{k+1} + A_k x_k}{A_{k +1}} - y_{k + 1}\right\rangle \\&+& \frac{L_{k + 1}}{2}\left\|\frac{\alpha_{k + 1}u_{k+1} + A_k x_k}{A_{k +1}} - \frac{\alpha_{k + 1}u_k + A_k x_k}{A_{k +1}}\right\|^2_{2} + 2\delta_k\\ &=& \frac{A_k }{A_{k +1}}\left(f_{\delta_k, L, \mu}(y_{k + 1}) + \left\langle g_{\delta,L, \mu}(y_{k + 1}), x_k - y_{k + 1}\right\rangle\right) \\&+&\frac{\alpha_{k + 1}}{A_{k +1}} \left(f_{\delta_k, L, \mu}(y_{k + 1}) + \left\langle g_{\delta,L, \mu}(y_{k + 1}), u_{k+1} - y_{k + 1}\right\rangle\right) + \frac{L_{k + 1}\alpha^2_{k + 1}}{2 A^2_{k +1}}\left\|u_{k+1} - u_k \right\|^2_{2} + 2\delta_k \\ &=& \frac{A_k }{A_{k +1}}\left(f_{\delta_k, L, \mu}(y_{k + 1}) + \left\langle g_{\delta,L, \mu}(y_{k + 1}), x_k - y_{k + 1}\right\rangle\right) \\&+&\frac{\alpha_{k + 1}}{A_{k +1}} \left(f_{\delta_k, L, \mu}(y_{k + 1}) + \left\langle g_{\delta,L, \mu}(y_{k + 1}), u_{k+1} - y_{k + 1}\right\rangle +\frac{1 + A_k\mu}{2\alpha_{k+1}}\left\|u_{k+1} - u_k \right\|^2_{2}\right) + 2\delta_k\\
    &\leq& \frac{A_k }{A_{k +1}}f_{\delta_k, L, \mu}(x_{k }) \\&+&\frac{\alpha_{k + 1}}{A_{k +1}} \left(f_{\delta_k, L, \mu}(y_{k + 1}) + \left\langle g_{\delta,L, \mu}(y_{k + 1}), u_{k+1} - y_{k + 1}\right\rangle +\frac{1 + A_k\mu}{2\alpha_{k+1}}\left\|u_{k+1} - u_k \right\|^2_{2}\right) + 2\delta_k\\
\end{eqnarray*}

Due to $u_{k+1}$ is the solution to $\argmin\limits_{x \in \mathcal{X}} \phi_{k + 1}(x)$, using \eqref{lemma stonyakin}  we have for any $x \in \mathcal{X}$ with probability at least $1-\sigma$:
\begin{eqnarray*}
    &&\alpha_{k + 1}\left\langle g_{\delta,L, \mu}(y_{k + 1}), u_{k+1} - y_{k + 1}\right\rangle + \frac{(1 +  A_k\mu)}{2}\|u_{k+1} - u_k\|^2_2\\ &+& \frac{\alpha_{k+1}\mu}{2}\|u_{k+1}  - y_{k + 1}\|^2_2 + \frac{(1 +  A_{k + 1}\mu)}{2}\|u_{k+1} - x\|^2_2 \\&\leq& 
    \alpha_{k + 1}\left\langle g_{\delta,L, \mu}(y_{k + 1}), x - y_{k + 1}\right\rangle + \frac{(1 +  A_k\mu)}{2}\|x - u_k\|^2_2 + \frac{\alpha_{k+1}\mu}{2}\|x  - y_{k + 1}\|^2_2 
\end{eqnarray*}
Then for all $x \in \mathcal{X}$ we have with probability at least $1-\sigma$:
\begin{eqnarray*}
    &&\alpha_{k + 1}\left\langle g_{\delta,L, \mu}(y_{k + 1}), u_{k+1} - y_{k + 1}\right\rangle + \frac{(1 +  A_k\mu)}{2}\|u_{k+1} - u_k\|^2_2 \\&\leq& 
    \alpha_{k + 1}\left\langle g_{\delta,L, \mu}(y_{k + 1}), x - y_{k + 1}\right\rangle + \frac{(1 +  A_k\mu)}{2}\|x - u_k\|^2_2 -  \frac{(1 +  A_{k + 1}\mu)}{2}\|u_{k+1} - x\|^2_2 \\&+& \frac{\alpha_{k+1}\mu}{2}\|x  - y_{k + 1}\|^2_2 
\end{eqnarray*}
Then for all $x \in \mathcal{X}$ we have with probability at least $1-\sigma$:
\begin{eqnarray*}
    f(x_{k + 1}) &\leq& \frac{A_k }{A_{k +1}}f(x_{k }) + \frac{\alpha_{k + 1}}{A_{k +1}} f(x) + 2\delta_k + \frac{\alpha_{k + 1}}{A_{k +1}} \Big( \frac{1 +  A_{k + 1}\mu}{2\alpha_{k+1}}\|u_{k+1} - x\|^2_2 +\frac{1 + A_k\mu}{2\alpha_{k+1}}\left\|x - u_k \right\|^2_{2} \Big) \\ &=& \frac{A_k }{A_{k +1}}f(x_{k }) + \frac{\alpha_{k + 1}}{A_{k +1}} f(x) + 2\delta_k + \frac{1 +  A_{k + 1}\mu}{2A_{k +1}}\|u_{k+1} - x\|^2_2 +\frac{1 + A_k\mu}{2A_{k +1}}\left\|x - u_k \right\|^2_{2} 
\end{eqnarray*}

Now we prove the statement of this theorem. First, we note that if we sum the inequality \eqref{convergence lemma} by $k$ from $0$ to $N-1$, for any $x \in \mathcal{X}$ we get 
\begin{eqnarray*}
    A_{N}f(x_{N}) - A_{N}f(x^*)  &\leq&  \sum^{N-1}_{k=0}2\delta_{k+1}A_{k+1} -\frac{1 +  A_{N}\mu}{2}\|x - u_{N}\|^2_2 + \frac{1}{2}\|x - u_0\|^2_2
\end{eqnarray*}
Then using the resulting inequality above for any $x \in \mathcal{X}$ we get
\begin{eqnarray*}
   \mathbb{P}\Bigg\{ A_{N}f(x_{N}) &-& A_{N}f(x^*) > \sum^{N-1}_{k=0}2\delta_{k+1}A_{k+1} -\frac{1 +  A_{N}\mu}{2}\|x - u_{N}\|^2_2 + \frac{1}{2}\|x - u_0\|^2_2\Bigg\}\\ &\leq& \sum^{N-1}_{k = 0}\mathbb{P}\Bigg\{ A_{k+1}f(x_{k +1}) - A_{k}f(x_{k}) + \frac{(1 +  A_{k+1}\mu)}{2}\|x - u_{k+1}\|^2_2 \\&&- \frac{(1 +  A_k\mu)}{2}\|x - u_k\|^2_2  \leq \alpha_{k+1} f(x) + 2\delta_{k+1}A_{k+1}\Bigg\} < N\sigma
\end{eqnarray*}

In this way, for any $x \in \mathcal{X}$ with  probability at least $(1 - N\sigma)$ we have 
\begin{eqnarray*}
    A_{N}f(x_{N}) - A_{N}f(x^*)  &\leq&  \sum^{N-1}_{k=0}2\delta_{k+1}A_{k+1} + \frac{1}{2}\|x - u_0\|^2_2
\end{eqnarray*}
We finish proof with next inequality (see \cite{stonyakin2020inexact}):
\begin{equation*}
    A_N \geq \frac{1}{2L}\exp\left(\frac{N - 1}{2}\sqrt{\frac{\mu}{L}}\right)
\end{equation*}
\end{proof}

\begin{corollary}
\label{fgm corollary}
Let function $f$ be convex on convex set $\mathcal{X}$ and be  equipped with a first-order $(\delta, \frac{\sigma}{N}, L, \mu)$-oracle. If the sequence $\{\delta\}_{k \geq 0}$ is bounded by $\delta$, we have with probability at least $(1 - \sigma)$:
\begin{eqnarray}
\label{fgm convergence 2}
    f(x_N) - f(x^*)&\leq&2L\exp\left(-\frac{N - 1}{2}\sqrt{\frac{\mu}{L}}\right) R^2_2 + \left(1 +\sqrt{\frac{L}{\mu}}\right)\delta
\end{eqnarray}
where $R_2$ is such that $\frac{1}{2}\|x_0 - x^*\|^2_2 \leq R^2_2$ and $x_0$ is the starting point.
\end{corollary}
\begin{proof}
To prove this statement, we give an auxiliary lemma
\begin{lemma} [see \cite{Devolder2013FirstorderMW}]
\label{devolder}
The sequence $\{A_k\}_{k \leq 0}$ satisfies
\begin{equation}
    \frac{\sum^k_{i = 0}A_i}{A_k} \leq 1 +\sqrt{\frac{L}{\mu}}
\end{equation}
\end{lemma}
We get the result of the Corollary immediately using \eqref{fgm convergence} and Lemma \ref{devolder}.
\end{proof}

\subsubsection{Analysis of Approach \ref{large scale minimax approach}}
For further analysis, we present the main lemma of this subsection
\begin{lmmrep}[see \cite{alkousa2020accelerated}]
\label{main lemma}
We denote \begin{equation*}
    y^*_f(x)=\argmax_{y \in \mathbb{R}^{n_y}}f(x, y),~  x^*_f(y)=\argmin_{x \in \mathcal{X}}f(x, y).
\end{equation*}Under assumption 1(s), 2 we have 
\begin{itemize}
    \item Function $x^*_f(y)$ is $\frac{L_{xy}}{\mu_x}$ -Lipschitz continuous w.r.t. the norm $\|\cdot\|_{2}$.
    \item Function $y^*_f(x)$ is $\frac{L_{xy}}{\mu_y}$ -Lipschitz continuous w.r.t. the norm $\|\cdot\|_{2}$
    \item Function $g(x)$ (see \eqref{inner}) is $L_g:=\left(L_{xx} + \frac{2L^2_{xy}}{\mu_y}\right)$-smooth w.r.t. the norm $\|\cdot\|_{2}$.
    \item Let $\widetilde{y}_{\delta}(x)$ be a $(\delta, \sigma)$-solution to problem $\max\limits_{y \in \mathbb{R}^{n_y}}f(x, y)$. Then,  for any $x', x \in \mathcal{X}$,  with probability at least  $1 - \sigma$ we have:
    \begin{equation*}
        \frac{\mu_x}{2}\|x - x'\|^2_{2} \leq g(x') - f(x,\widetilde{y}_{\delta}(x)) -\left\langle \nabla_x f(x, \widetilde{y}_{\delta}(x)),x' -x \right\rangle \leq \frac{2L_g}{2}\|x' - x\|^2_{2} + 2\delta.
    \end{equation*}
    \item We define
    \begin{equation*}
        g(x)=\max_{y \in \mathbb{R}^{n_y}} f(x, y),~h(y)=\min_{x \in \mathcal{X}} f(x, y)
    \end{equation*}
    Let $\hat{x}$ be $(\varepsilon_x, \sigma_x)$-solution of  problem  $\min\limits_{x \in \mathcal{X}}g(x)$, let $\widetilde{y}_{\varepsilon_y}(\hat{x})$ be $(\varepsilon_y, \sigma_y)$-solution of  problem $\max\limits_{y \in \mathbb{R}^{n_y}}f(\hat{x}, y)$. Then $\widetilde{y}_{\varepsilon_y}(\hat{x})$ is $(\tilde{\varepsilon}, 1-\sigma_x - \sigma_y)$-solution to problem \eqref{problem}, where
    \begin{equation*}
        \tilde{\varepsilon} = \left(\frac{L_{yy}}{\mu_y} + \frac{2L^2_{xy}}{\mu_x\mu_y}\right) \varepsilon_y + \left(\frac{L^2_{xy}L_{yy}}{\mu_x\mu^2_y} + \frac{2L^4_{xy}}{\mu^2_x\mu^2_y} \right) \varepsilon_x.
    \end{equation*}
\end{itemize}
\end{lmmrep}
\begin{proof}
To justify the statement of the lemma, we give a more general statement.
\begin{lemma}[see lemma \ref{main lemma} ]
We denote \begin{equation*}
    y^*_f(x)=\argmax_{y \in \mathbb{R}^{n_y}}g(x, y);~  x^*_f(y)=\argmin_{x \in \mathcal{X}}g(x, y).
\end{equation*} Under assumption 1(s) for $\|\cdot\|_{p_x}$ and $\|\cdot\|_{p_y}$ norms, assumption 2 for $\|\cdot\|_{p_x}$ and $\|\cdot\|_{p_y}$ norms we have 
\begin{itemize}
    \item Function $x^*_g(y)$ is $\frac{L_{xy}}{\mu_x}$ -Lipschitz continuous w.r.t. the norm $\|\cdot\|_{p_x}$.
    \item Function $y^*_g(x)$ is $\frac{L_{xy}}{\mu_y}$ -Lipschitz continuous w.r.t. the norm $\|\cdot\|_{p_x}$
    \item Function $g(x)$ is $L_g=\left(L_{xx} + \frac{2L^2_{xy}}{\mu_y}\right)$-smooth w.r.t. the norm $\|\cdot\|_{p_x}$.
    \item Let $\widetilde{y}_{\delta}(x)$ - $(\delta, \sigma)$ be solution to problem $\max\limits_{y \in \mathbb{R}^{n_y}}f(x, y)$. Then  for any $x', x \in \mathcal{X}$  with probability at least  $1 - \sigma$ we have:
    \begin{equation*}
        \frac{\mu_x}{2}\|x - x'\|^2_{p_x} \leq g(x') - f(x,\widetilde{y}_{\delta}(x)) -\left\langle \nabla_x f(x, \widetilde{y}_{\delta}(x)),x' -x \right\rangle \leq \frac{2L_g}{2}\|x' - x\|^2_{p_x} + 2\delta
    \end{equation*}
    \item We define
    \begin{equation*}
        g(x)=\max_{y \in \mathbb{R}^{n_y}} f(x, y),~h(y)=\min_{x \in \mathcal{X}} f(x, y)
    \end{equation*}
    Let $\hat{x}$ be $(\varepsilon_x, \sigma_x)$-solution of  problem  $\min\limits_{x \in \mathcal{X}}g(x)$, let $\widetilde{y}_{\varepsilon_y}(\hat{x})$ be $(\varepsilon_y, \sigma_y)$-solution of  problem $\max\limits_{y \in \mathbb{R}^{n_y}}f(\hat{x}, y)$. Then $\widetilde{y}_{\varepsilon_y}(\hat{x})$ is $(\tilde{\varepsilon}, 1-\sigma_x - \sigma_y)$-solution to problem, where
    \begin{equation*}
        \tilde{\varepsilon} = \left(\frac{L_{yy}}{\mu_y} + \frac{2L^2_{xy}}{\mu_x\mu_y}\right) \varepsilon_y + \left(\frac{L^2_{xy}L_{yy}}{\mu_x\mu^2_y} + \frac{2L^4_{xy}}{\mu^2_x\mu^2_y} \right) \varepsilon_x
    \end{equation*}
\end{itemize}
\end{lemma}
\textbf{Proof}\begin{enumerate}
    \item Estimate $\| y^*(x) - y^*(x')\|_{p_y}$:
    \begin{eqnarray*}
        \| y^*(x) &-& y^*(x')\|^2_{p_y} \leq \frac{2}{\mu_y}(f(x, y^*(x)) - f(x, y^*(x')))\\&\leq& \frac{2}{\mu_y}(f(x, y^*(x)) - f(x, y^*(x'))) - \frac{2}{\mu_y}(f(x', y^*(x)) - f(x', y^*(x')))\\ &=&\frac{2}{\mu_y}\int\limits^{1}_{0}\langle \nabla_x f(x' +t(x - x'), y^*(x)) - \nabla_x f(x' +t(x- x'), y^*(x')), x' -x\rangle dt \\ &\leq& \frac{2}{\mu_y}\int\limits^{1}_{0}\|\nabla_x f(x' +t(x - x'), y^*(x)) - \nabla_x f(x' +t(x- x'), y^*(x'))\|_{q_x}\| x' -x\|_{p_x} dt \\ &\leq& \frac{2}{\mu_y} L_{xy} \| x' -x\|_{p_x} \| y^*(x) - y^*(x')\|_{p_y}
    \end{eqnarray*}
    In thus way, we have
    \begin{equation*}
        \| y^*(x) - y^*(x')\|_{p_y}  \leq \frac{2L_{xy}}{\mu_y}\|x'-x\|_{p_x} 
    \end{equation*}
    
    \item Estimate $\| x^*(y) - x^*(y')\|_{p_x}$. Due to $\mu_x$-strong convexity function $f(x, y)$:
    \begin{eqnarray*}
        \| x^*(y) &-& x^*(y')\|^2_{p_x} \leq \frac{2}{\mu_x}(f(x^*(y), y') - f(x^*(y'), y'))\\&\leq& \frac{2}{\mu_x}(f(x^*(y), y') - f(x^*(y'), y')) - \frac{2}{\mu_x}(f(x^*(y), y) - f(x^*(y'), y))\\ &=&\frac{2}{\mu_x}\int\limits^{1}_{0}\langle \nabla_y f(x^*(y),y +t(y' - y)) - \nabla_y f( x^*(y'), y +t(y'- y)), y' -y\rangle dt \\ &\leq& \frac{2}{\mu_x}\int\limits^{1}_{0}\|\nabla_y f(x^*(y),y +t(y' - y)) - \nabla_y f( x^*(y'), y +t(y'- y))\|_{q_y}\| y' -y\|_{p_y} dt \\ &\leq& \frac{2}{\mu_x} L_{xy} \| y' -y\|_{p_y} \| x^*(y) - x^*(y')\|_{p_x}
    \end{eqnarray*}
    We have
    \begin{equation*}
        \| x^*(y) - x^*(y')\|_{p_x}  \leq \frac{2L_{xy}}{\mu_x}\|y'-y\|_{p_y} 
    \end{equation*}
    
    \item Consider the following value $\|\nabla g(x) - \nabla g(x')\|_{q_x}$:
    \begin{eqnarray*}
        \|\nabla g(x) - \nabla g(x')\|_{q_x} &=& \|\nabla_x f(x, y^*(x)) - \nabla_x f(x', y^*(x'))\|_{q_x} \\&\leq& \|\nabla_x f(x, y^*(x)) - \nabla_x f(x, y^*(x'))\|_{q_x}\\ &+& \|\nabla_x f(x, y^*(x')) - \nabla_x f(x', y^*(x'))\|_{q_x} \\&\leq& L_{xx} \|x - x'\|_{p_x} + L_{xy}\| y^*(x) - y^*(x')\|_{p_y}
    \end{eqnarray*}
    Then we have:
    \begin{eqnarray*}
        \|\nabla g(x) - \nabla g(x')\|_{q_x} &\leq& \left(L_{xx} + \frac{2L^2_{xy}}{\mu_y}\right)\|x - x'\|_{p_x}
    \end{eqnarray*}
    
    \item  Let $\widetilde{y}_{\delta}(x)$ - $(\delta, \sigma)$ be solution to problem $\max\limits_{y\in \mathbb{R}^{n_y}} f(x, y),$ and $y^*(x)$ be exact solution. We prove the left side of  inequality. For any $x, x' \in \mathcal{X}$ we have
    \begin{eqnarray*}
        g(x') = f(x', y^*(x')) \geq f(x', \widetilde{y}_{\delta}(x)) \geq f(x, \widetilde{y}_{\delta}(x)) + \left\langle \nabla_x f(x, \widetilde{y}_{\delta}(x)),x' -x  \right\rangle +\frac{\mu_x}{2}\|x - x'\|^2_{p_x}
    \end{eqnarray*}
    Then we have
    \begin{equation*}
        \frac{\mu_x}{2}\|x - x'\|^2_{p_x} \leq g(x') - f(x,\widetilde{y}_{\delta}(x)) -\left\langle \nabla_x f(x, \widetilde{y}_{\delta}(x)),x' -x \right\rangle
    \end{equation*}
    Now we prove the left side of inequality. For any $x, x' \in \mathcal{X}$ we have with probability at least $1 - \sigma$ 
    \begin{eqnarray*}
        g(x') &\leq& g(x) +\left\langle\nabla g(x), x' - x \right\rangle  + \frac{L_{g}}{2}\|x - x'\|^2_{p_x} \\ &\leq& f(x,\widetilde{y}_{\delta}(x)) +  \left\langle\nabla_x f(x,\widetilde{y}_{\delta}(x)) , x' - x \right\rangle + \frac{L_g}{2}\|x - x'\|^2_{p_x} + \delta \\&+&\left\langle\nabla_x f(x,\widetilde{y}_{\delta}(x)) - \nabla g(x) , x - x' \right\rangle  \\ &\leq& f(x,\widetilde{y}_{\delta}(x)) +  \left\langle\nabla_x f(x,\widetilde{y}_{\delta}(x)) , x' - x \right\rangle + \frac{L_g}{2}\|x - x'\|^2_{p_x} + \delta \\&+&\left\|\nabla_x f(x,\widetilde{y}_{\delta}(x)) - \nabla g(x)\right\|_{q_x} \|x - x'\|_{p_x}
    \end{eqnarray*}
    Estimate $\left\|\nabla_x f(x,\widetilde{y}_{\delta}(x)) - \nabla g(x)\right\|_{q_x}$. we have with probability at least $1 - \sigma$ 
    \begin{eqnarray*}
        \left\|\nabla_x f(x,\widetilde{y}_{\delta}(x)) - \nabla g(x)\right\|^2_{q_x} &=& \left\|\nabla_x f(x,\widetilde{y}_{\delta}(x)) - \nabla_x f(x, y^*(x))\right\|^2_{q_x} \leq L^2_{xy}\|\widetilde{y}_{\delta}(x) - y^*(x)\|^2_{p_y}\\
        &\leq& \frac{2L^2_{xy}}{\mu_y} \left(f(x, y^*(x)) - f(x, \widetilde{y}_{\delta}(x))\right) \leq \frac{2L^2_{xy}}{\mu_y} \delta
    \end{eqnarray*}
    Then for any $x, x' \in \mathcal{X}$ we have with probability at least $1 - \sigma$  
    \begin{eqnarray*}
        g(x') &\leq& f(x,\widetilde{y}_{\delta}(x)) +  \left\langle\nabla_x f(x,\widetilde{y}_{\delta}(x)) , x' - x \right\rangle + \frac{L_g}{2}\|x - x'\|^2_{p_x} + \delta + \sqrt{\frac{2L^2_{xy}}{\mu_y} \delta}\|x - x'\|_{p_x}\\
        &\leq& f(x,\widetilde{y}_{\delta}(x)) +  \left\langle\nabla_x f(x,\widetilde{y}_{\delta}(x)) , x' - x \right\rangle + \left(\frac{L_g}{2} + \frac{2L^2_{xy}}{\mu_y}\right)\|x - x'\|^2_{p_x} + 2\delta \\ &\leq& f(x,\widetilde{y}_{\delta}(x)) +  \left\langle\nabla_x f(x,\widetilde{y}_{\delta}(x)) , x' - x \right\rangle + \left(\frac{L_{xx}}{2} + \frac{3L^2_{xy}}{\mu_y}\right)\|x - x'\|^2_{p_x} + 2\delta\\ &\leq& f(x,\widetilde{y}_{\delta}(x)) +  \left\langle\nabla_x f(x,\widetilde{y}_{\delta}(x)) , x' - x \right\rangle + 2\left(L_{xx} + \frac{2L^2_{xy}}{\mu_y}\right)\|x - x'\|^2_{p_x} + 2\delta
    \end{eqnarray*}
    I this way, we have function $g(x)$ equips  $(2\delta, \sigma,  2L_g, \mu_x)$-oracle.
    \item According to $\mu_x$-strong convexity of function $g(x)$,  we have with probability at least $1-\sigma_x$ 
    $$\|\hat{x} - x^*\|^2_{p_x} \leq \frac{2}{\mu_x} \varepsilon_x$$
    According to $\mu_x$-strong concavity of function $h(x)$, using union bound we have with probability at least $1-\sigma_y - \sigma_x$
    \begin{eqnarray*}
        \|\widetilde{y}_{\varepsilon_y}(\hat{x}) - y^*\|^2_{p_y} &\leq& 2\|\widetilde{y}_{\varepsilon_y}(\hat{x}) - y^*(\hat{x})\|^2_{p_y} + 2\|y^*(\hat{x}) - y^*\|^2_{p_y} \\&\leq& \frac{2}{\mu_y} \varepsilon_y + \frac{L^2_{xy}}{\mu^2_y}\|\hat{x} - x^*\|^2_{p_x} \leq \frac{2}{\mu_y} \varepsilon_y + \frac{2L^2_{xy}}{\mu_x\mu^2_y}\varepsilon_x
    \end{eqnarray*} 
    We have with probability at least $1-\sigma_y - \sigma_x$
    \begin{eqnarray*}
        h(y^*) - h\left(\widetilde{y}_{\varepsilon_y}(\hat{x})\right) &\leq& \frac{L_{yy} + \frac{2L^2_{xy}}{\mu_x}}{2}\|y^* - \widetilde{y}_{\varepsilon_y}(\hat{x})\|^2_{p_y} \\&\leq& \left(\frac{L_{yy}}{\mu_y} + \frac{2L^2_{xy}}{\mu_x\mu_y}\right) \varepsilon_y + \left(\frac{L^2_{xy}L_{yy}}{\mu_x\mu^2_y} + \frac{2L^4_{xy}}{\mu^2_x\mu^2_y} \right) \varepsilon_x
\end{eqnarray*}
\end{enumerate}
\end{proof}

Now we are ready to present the main result of this section
\begin{thmrep}
Let $\varepsilon >0$ be the target accuracy of the solution to the problem \eqref{problem} and $\sigma \in (0, 1)$ be the target confidence level. Let the auxiliary problems \eqref{inner}, \eqref{outer} be solved with accuracies 
\begin{equation*}
    \varepsilon_x = \widetilde{O}\left(\varepsilon\left(\frac{L^2_{xy} L_{yy}}{\mu_x(\mu_y + L_{yy})^2} + \frac{2L^4_{xy}}{\mu^2_x(\mu_y + L_{yy})^2} \right)^{-1}\right);
\end{equation*}
\begin{equation*}
    \varepsilon_y = \widetilde{O}\left(\varepsilon\left(\frac{L_{yy}}{\mu_y + L_{yy}} + \frac{2L^2_{xy}}{\mu_x(\mu_y + L_{yy})}\right)^{-1}\left(\frac{L_{xx}}{\mu_x} + \frac{2L^2_{xy}}{\mu_x(\mu_y + L_{yy})}\right)^{-\nicefrac{1}{2}}\right)
\end{equation*}
and confidence levels
\begin{equation*}
    \sigma_x = \widetilde{O}\left(\sigma\sqrt{\frac{\mu_y}{L_{yy}}}\right);
\end{equation*}
\begin{equation*}
    \sigma_y = \widetilde{O}\left(\sigma \left(\frac{L_{xx}L_{yy}}{\mu_x \mu_y} + \frac{2L^2_{xy} }{\mu_x \mu_y}\right)^{-\nicefrac{1}{2}}\right),
\end{equation*}
that is, a $(\varepsilon_x, \sigma_x)-$solution to the problem \eqref{outer} and a $(\varepsilon_y, \sigma_y)-$solution to the problem \eqref{inner} are found (see Definition \ref{es solution}).
Then, under assumptions 1(s), 2, the proposed Approach \ref{large scale minimax approach} guarantees to find an $(\varepsilon, \sigma)$-solution to the problem \eqref{problem}.
Moreover, the required number of calls to the first-order oracle $\nabla_x f(x, y)$ and  the zeroth-order oracle $f(x, y)$ satisfy the following bounds
    \begin{equation*}
        \text{Total  Number of Calls for $\nabla_x f(x, y)$ is }~ \widetilde{O}\left(\sqrt{\frac{L_{xx}}{\mu_x} +\frac{2L^2_{xy}}{\mu_x\mu_y  }}\right),
    \end{equation*}
    \begin{equation*}
        \text{Total  Number of Calls for $ f(x, y)$ is }~ \widetilde{O}\left(n_y\sqrt{\frac{L_{xx} L_{yy}}{\mu_x\mu_y} +\frac{2L^2_{xy}}{\mu_x\mu_y} }  \right).
    \end{equation*}
\end{thmrep}
\begin{proof}

{\bf Oracles  Complexities Analysis}

Now we can analyze the oracle complexity of our approach
\begin{enumerate}
    \item The number of iterations of the Catalyst algorithm \ref{Catalyst} (see \cite{catalyst}) to solve the problem \eqref{outer catalyst} with the accuracy of $\varepsilon$:
    \begin{equation}
    \label{catalyst iteration}
        \widetilde{O}\left(\sqrt{\frac{H_1}{\mu_y }}\right)
    \end{equation}
    
    \item Now we need to determine the number of iterations of the $\mathcal{M}$ method to solve the problem \eqref{subpr}. According to theorem \ref{ARDDsc convergence}, the number of iterations of algorithm \ref{ARDDsc} required to find the solution with accuracy $\varepsilon_y$ is 
    \begin{equation}
    \label{iterations of inner problem}
        N_{\text{ARDDsc}} = \widetilde{O}\left(n_{y}\sqrt{\frac{L_{yy} + H_1}{\mu_y + H_1 }} \log\left(\frac{(\mu_y + H_1)R^2_{p_y}}{\varepsilon_y}\right)\right).
    \end{equation}
    
    \item Now we solve the outer problem \eqref{catalyst outer} with the algorithm \ref{Fast adaptive gradient method}. We need to determine the Lipschitz constant and the strong convexity constant. The strong convexity constant remains the same, the Lipschitz constant, according to Lemma \ref{main lemma}, is $$L_g = L_{xx} + \frac{2L^2_{xy}}{\mu_y + H_1}.$$ Complexity of solution of the outer problem \eqref{catalyst outer} (number of calculations of the gradient value $\nabla_x f(x, y)$) with accuracy $\varepsilon_x$: 
    \begin{equation}
    \label{iterations of outer problem}
        N_{\text{outer}} = \widetilde{O}\left(\sqrt{\frac{L_{g}}{\mu_x}} \right) = \widetilde{O}\left(\sqrt{\frac{L_{xx}}{\mu_x} +\frac{2L^2_{xy}}{\mu_x(\mu_y + H_1)} }\right)
    \end{equation}
    
    \item Thus, the total complexity of calls of zeroth-order oracle $f(x, y)$, solving the inner problem \eqref{problem} by the algorithm \ref{ARDDsc} with accuracy $\varepsilon$, taking $H_1 = L_{yy}$, is equal to
    \begin{eqnarray*}
        N_{\text{inner}} &=& \widetilde{O}\left(n_y\sqrt{\frac{L_{yy}+ H_1}{\mu_y + H_1 }}\sqrt{\frac{H_1}{\mu_y}}\sqrt{\frac{L_{xx}}{\mu_x} +\frac{2L^2_{xy}}{\mu_x(\mu_y + H_1)} }  \right)\\ &=& \widetilde{O}\left(n_y\sqrt{\frac{L_{yy}+ L_{yy}}{\mu_y + L_{yy} }}\sqrt{\frac{L_{yy}}{\mu_y}}\sqrt{\frac{L_{xx}}{\mu_x} +\frac{2L^2_{xy}}{\mu_x(\mu_y + L_{yy})} }  \right)\\ &=& \widetilde{O}\left(n_y\sqrt{\frac{L_{yy}+ L_{yy}}{ L_{yy} }}\sqrt{\frac{L_{xx} L_{yy}}{\mu_x\mu_y} +\frac{2L^2_{xy}L_{yy}}{\mu_x\mu_y L_{yy}} }  \right)\\ &=& \widetilde{O}\left(n_y\sqrt{\frac{L_{xx} L_{yy}}{\mu_x\mu_y} +\frac{2L^2_{xy}}{\mu_x\mu_y} }  \right) 
    \end{eqnarray*} 
\end{enumerate}

{\bf Inexactness  Complexities Analysis}

Now it remains to determine with what accuracy $\varepsilon_x$, $\varepsilon_y$ it is necessary to solve inner \eqref{catalyst inner} and outer \eqref{catalyst outer} problems, so that the number of calls to the oracles satisfies this ratio \eqref{Oracle Complexity}. In order to answer this question, we will introduce new designations and redefine some things. First, we denote
\begin{equation*}
    \xi(x)=\max_{y \in \mathbb{R}^{n_y}}\psi(x, y),~  \zeta(y)=\min_{x \in \mathcal{X}}\psi(x, y)
\end{equation*}
Let $\hat{x} \in \mathcal{X}$ be $(\varepsilon_x, \sigma_x)$-solution to problem $\min\limits_{x \in \mathcal{X}}\xi(x)$, in other words,  we have with probability at least $1 -\sigma_x$:
\begin{equation*}
    0 \leq \xi(\hat{x}) - \min\limits_{x \in \mathcal{X}}\xi(x) = \psi(\hat{x}, y^*(\hat{x})) -   \psi(x^*, y^*) \leq \varepsilon_x
\end{equation*}
Let $\widetilde{y}_{\varepsilon_y}(\hat{x})$ - $(\varepsilon_y, \sigma_y)$-solution to problem $\max\limits_{y \in \mathbb{R}^{n_y}}\psi(\hat{x}, y)$, in other words,  we have with probability at least $1 -\sigma_y$:
\begin{equation*}
     0 \leq \psi(\hat{x}, y^*(\hat{x})) -   \psi(\hat{x}, \widetilde{y}_{\varepsilon_y}(\hat{x})) \leq \varepsilon_y
\end{equation*}

First, according to corollary \ref{fgm corollary} for $\min\limits_{x \in \mathcal{X}}\xi(x)$, we have a limit on noise $\varepsilon_y$:
\begin{equation*}
     \varepsilon_y\left(1+ \sqrt{\frac{L_g}{\mu_x}} \right) = O\left(\varepsilon_{x}\right)
\end{equation*}
Thus, we have 
\begin{equation*}
    \varepsilon_y = O\left(\varepsilon_x\sqrt{\frac{\mu_x}{L_g}}\right)
\end{equation*}
It remains only to understand how to define $\varepsilon_{x}$ via $\varepsilon$.To answer this question using the main lemma \ref{main lemma}, replacing $\mu_y$ with $\mu_y+H_1$ and $L_{yy}$ with $L_{yy} + H_1$, we obtain inequality:
\begin{equation*}
    \tilde{\varepsilon}\geq  \left(\frac{L_{yy} + H_1}{\mu_y + H_1} + \frac{2L^2_{xy}}{\mu_x\mu_y}\right) \varepsilon_y + \left(\frac{L^2_{xy}(L_{yy} + H_1)}{\mu_x(\mu_y + H_1)^2} + \frac{2L^4_{xy}}{\mu^2_x(\mu_y + H_1)^2} \right) \varepsilon_x
\end{equation*}
According to results from \cite{ivanova2020adaptive}, we have  $\tilde{\varepsilon} = O(\varepsilon)$. Then we have
\begin{eqnarray*}
    \varepsilon_x & = & \widetilde{O}\left(\varepsilon\left(\frac{L^2_{xy}(L_{yy}+H_1)}{\mu_x(\mu_y + H_1)^2} + \frac{2L^4_{xy}}{\mu^2_x(\mu_y + H_1)^2} \right)^{-1}\right)
    \\ \varepsilon_y & = & \widetilde{O}\left(\varepsilon\left(\frac{L_{yy} + H_1}{\mu_y + H_1} + \frac{2L^2_{xy}}{\mu_x(\mu_y + H_1)}\right)^{-1}\sqrt{\frac{\mu_x}{L_g}}\right)
\end{eqnarray*}
Using that $H_1 = L_{yy}$, we have 
\begin{eqnarray*}
    \varepsilon_x & = & \widetilde{O}\left(\varepsilon\left(\frac{L^2_{xy} L_{yy}}{\mu_x(\mu_y + L_{yy})^2} + \frac{2L^4_{xy}}{\mu^2_x(\mu_y + L_{yy})^2} \right)^{-1}\right)
    \\ \varepsilon_y & = & \widetilde{O}\left(\varepsilon\left(\frac{L_{yy}}{\mu_y + L_{yy}} + \frac{2L^2_{xy}}{\mu_x(\mu_y + L_{yy})}\right)^{-1}\sqrt{\frac{\mu_x}{L_{xx} + \frac{2L^2_{xy}}{\mu_y + L_{yy}}}}\right)
\end{eqnarray*}

According to theorem \ref{fgm theorem} and approach \ref{large scale minimax approach}, we have that 
\begin{equation*}
    \sigma_x  = N_{\text{outer}} \sigma_y, \text {where } N_{\text{outer}} = \widetilde{O}\left(\sqrt{\frac{L_{xx}}{\mu_x} + \frac{2L^2_{xy}}{\mu_x (\mu_y + L_{yy}) }}\right)
\end{equation*} 
Having obtained this ratio, we now need to determine how they relate to each other. Using the catalyst algorithm and the convergence proof, one can easily show that 
\begin{equation*}
    \sigma  \geq N_{\text{Catalyst}} (\sigma_y + \sigma_x) > N_{\text{Catalyst}} N_{\text{outer}} \sigma_y  , \text{ where } N_{\text{Catalyst}} = \widetilde{O}\left(\sqrt{\frac{L_{yy}}{\mu_y}}\right)
\end{equation*} 
Thus we have 
\begin{eqnarray*}
    \sigma_y &=& \widetilde{O}\left(\frac{\sigma}{N_{\text{Catalyst}} N_{\text{outer}}}\right) =  \widetilde{O}\left(\sigma\sqrt{\frac{\mu_y}{L_{yy}} \left(\frac{L_{xx}}{\mu_x} + \frac{2L^2_{xy}}{\mu_x (\mu_y + L_{yy})}\right)^{-1}}\right) \\&=& \widetilde{O}\left(\sigma\sqrt{ \left(\frac{L_{xx}L_{yy}}{\mu_x \mu_y} + \frac{2L^2_{xy} L_{yy}}{\mu_x (\mu_y + L_{yy}) \mu_y}\right)^{-1}}\right) \\&=& \widetilde{O}\left(\sigma \left(\frac{L_{xx}L_{yy}}{\mu_x \mu_y} + \frac{2L^2_{xy} }{\mu_x \mu_y}\right)^{-\nicefrac{1}{2}}\right)
\end{eqnarray*}
\begin{eqnarray*}
    \sigma_x &=& \sigma_y N_{\text{outer}} =  \widetilde{O}\left(\frac{\sigma N_{\text{outer}}}{N_{\text{Catalyst}} N_{\text{outer}}}\right) = \widetilde{O}\left(\frac{\sigma }{N_{\text{Catalyst}}}\right) = \widetilde{O}\left(\sigma\sqrt{\frac{\mu_y}{L_{yy}}}\right)
\end{eqnarray*}
\end{proof}

%
%
%
\bibliographystyle{splncs04}
\bibliography{literature}

\appendix

\end{document}